\documentclass[a4paper]{amsart}

\title{On the $K$-theory of Higher Rank Graph $C^*$-Algebras}

\author{D. Gwion Evans}

\address{Institute of Mathematical and Physical Sciences\\ Aberystwyth University\\ Penglais Campus \\ Aberystwyth\\ Ceredigion \\ SY23 3BZ \\ Wales \\ UK.}
\email{dfe@aber.ac.uk}

\thanks{Research supported by the European Union Research Training Network in Quantum Spaces - Noncommutative Geometry.}



\usepackage{amssymb,amscd}
\newcommand{\field}[1]{\mathbb{#1}}

\newcommand{\T}{\field{T}}
\newcommand{\N}{\field{N}}
\newcommand{\Z}{\field{Z}}
\newcommand{\K}{\field{K}}
\newcommand{\KK}{\mathcal{K}}

\newcommand{\A}{\mathcal{A}}
\newcommand{\E}{\mathcal{E}}
\newcommand{\mbf}[1]{\mbox{\boldmath$#1$}}

\newcommand{\Span}{\overline{\spn}}
\newcommand{\Homo}{H}
\newcommand{\into}{\longrightarrow}
\newcommand{\outof}{\longleftarrow}
\DeclareMathOperator{\im}{im}
\DeclareMathOperator{\spn}{span}
\DeclareMathOperator{\coker}{coker}
\DeclareMathOperator{\tor}{tor}
\DeclareMathOperator{\aut}{\mbox{Aut}}

\newtheorem{thrm}{Theorem}[section]
\newtheorem{prop}[thrm]{Proposition}
\newtheorem{lemma}[thrm]{Lemma}
\newtheorem{cor}[thrm]{Corollary}
\theoremstyle{definition}
\newtheorem{defn}[thrm]{Definition}
\newtheorem{rem}[thrm]{Remarks}
\newtheorem{exs}[thrm]{Examples}
\newtheorem{convention}[thrm]{Convention}

\begin{document}
\maketitle
{\em To appear in the New York Journal of Mathematics (http://nyjm.albany.edu:8000/).}
\begin{abstract}
Given a row-finite $k$-graph $\Lambda$ with no sources we
investigate the $K$-theory of the higher rank graph $C^*$-algebra, $C^*(\Lambda)$. When $k=2$ we are able to give explicit
formulae to calculate the $K$-groups of $C^*(\Lambda)$. The $K$-groups of $C^*(\Lambda)$ for $k>2$
can be calculated under
certain circumstances and we consider the case $k=3$. We prove that for arbitrary $k$, the
torsion-free rank of $K_0(C^*(\Lambda))$ and $K_1(C^*(\Lambda))$
are equal when $C^*(\Lambda)$ is unital, and for $k=2$ we determine the position of
the class of the unit of $C^*(\Lambda)$ in $K_0(C^*(\Lambda))$.
\end{abstract}

\section{Introduction}

In \cite{S91} Spielberg realised that a crossed product
algebra $C(\Omega)\rtimes\Gamma$, where $\Omega$ is the boundary of
a certain tree and $\Gamma$ is a free group, is isomorphic to a
Cuntz-Krieger algebra \cite{CK80,C81}. Noticing that such a tree may be regarded as an
affine building of type $\tilde{A}_1$, Robertson and
Steger studied the situation when a group $\Gamma$ acts simply
transitively on the vertices of an affine building of type
$\tilde{A}_2$ with boundary $\Omega$ \cite{RS96}. They found that the
corresponding crossed product algebra $C(\Omega)\rtimes \Gamma$ is
generated by two Cuntz-Krieger algebras. This led them
to define a $C^*$-algebra $\A$ via a finite sequence of finite
0--1 matrices (i.e. matrices with entries in $\{0,1\}$)
$M_1,\ldots,M_r$ satisfying certain conditions (H0)-(H3), such that $\A$ is generated by $r$ Cuntz-Krieger
algebras, one for each $M_1,\ldots,M_r$. Accordingly they named
their algebras higher rank Cuntz-Krieger algebras, the rank being
$r$.

Kumjian and Pask \cite{KP00} noticed that Robertson and Steger had
constructed their algebras from a set, $W$ of {\it (higher rank)
words} in a finite {\it alphabet} $A$ - the common index set of
the 0--1 matrices - and realised that $W$ could be thought of
as a special case of a generalised directed graph - a higher rank graph.
Subsequently, Kumjian and Pask associated a $C^*$-algebra,
$C^*(\Lambda)$ to the higher rank graph $\Lambda$ and showed that
$\A\cong C^*(W)$ \cite[Corollary 3.5 (ii)]{KP00}. Moreover, they derived a number of results elucidating the structure
of higher rank graph $C^*$-algebras. They show in \cite[Theorem
5.5]{KP00} that a simple, purely infinite $k$-graph $C^*$-algebra $C^*(\Lambda)$ may be
classified by its $K$-theory. This is a consequence of $C^*(\Lambda)$
satsifying the hypotheses of the Kirchberg-Phillips
classification theorem (\cite{K,P00}). Moreover, criteria on the underlying
$k$-graph $\Lambda$ were found that decided when $C^*(\Lambda)$ was
simple and purely infinite (see \cite[Proposition 4.8, Proposition 4.9]{KP00} and \cite{S06}). Thus
a step towards the classification of $k$-graph $C^*$-algebras is the
computation of their $K$-groups.

In \cite[Proposition 4.1]{RS01} Robertson and Steger proved that the
$K$-groups of a rank 2 Cuntz-Krieger algebra is given in terms of the
homology of a certain chain complex, whoose differentials are defined
in terms of $M_1,\dots,M_r$. Their proof relied on the
fact that a rank 2 Cuntz-Krieger algebra is stably isomorphic to a
crossed product of an AF-algebra by $\Z^2$. We will
generalise their method to provide explicit formulae for the $K$-groups of 2-graph
$C^*$-algebras and to gain information on the $K$-groups of
$k$-graph $C^*$-algebras for $k>2$.

The rest of this paper is organised as follows. We begin in \S\ref{S:Prelim} by recalling
the fundamental definitions relating to higher rank graphs and their $C^*$-algebras we will need from \cite{KP00}.

In \S\ref{S:K-theory} we
use the fact that the $C^*$-algebra of a row-finite $k$-graph $\Lambda$ with no sources
is stably isomorphic to a crossed product of an AF algebra, $B$, by
$\Z^k$ (\cite[Theorem 5.5]{KP00}) to apply a theorem of Kasparov \cite[6.10 Theorem]{K88} to deduce that there is a homological spectral sequence
(\cite[Chapter 5]{W94})
converging to $K_*(C^*(\Lambda))$ with initial term given by
$E^2_{p,q}\cong H_p(\Z^k,K_q(B))$ (see \cite{Br94} for the definition
of the homology of a group $G$
with coefficients in a left $G$-module $M$, denoted by $H_*(G,M)$). We
will see that it
suffices to compute $H_*(\Z^k,K_0(B))$. It transpires that
$H_*(\Z^k,K_0(B))$ is given by the so called vertex matrices of
$\Lambda$. These are matrices with non-negative integer entries
that encode the structure of the category $\Lambda$. Next we assemble
the results of \S\ref{S:K-theory} and state them in our main theorem, Theorem
\ref{T:Main}. We then specialise to the cases $k=2$ and $k=3$. For
$k=2$ a complete description of the $K$-groups in terms of the vertex
matrices can be given. For $k=3$ we illustrate how Theorem
\ref{T:Main} can be used to give a description of the $K$-groups of
3-graph $C^*$-algebras under stronger hypotheses.

In section \S\ref{S:unital} we consider the $K$-theory of unital
$k$-graph $C^*$-algebras. We show that the torsion-free rank of
$K_0(C^*(\Lambda))$ is equal to that of $K_1(C^*(\Lambda))$ when
$C^*(\Lambda)$ is unital and give formulae for the torsion-free rank
and torsion parts of the $K$-groups of 2-graph
$C^*$-algebras.

We conclude with \S\ref{S:examples}, in which we consider some immediate applications to the classification of $k$-graph $C^*$-algebras by means of the Kirchberg-Phillips classification theorem.  We also consider some simple examples of $K$-group calculations using the results derived in the previous sections.

This paper was written while the author was an European Union Network in Quantum
Spaces -- Non-Commutative Geometry funded post-doc at the
University of Copenhagen. The paper develops a part of the author's PhD
thesis, which was written under the supervision of David
E. Evans at Cardiff University. We would like to take this opportunity
to thank David for his guidance and support, and Johannes
Kellendonk and Ryszard Nest for enlightening
discussions on homological algebra. We would also like to express our gratitude to the members of the operator
algebras groups in both universities, for maintaining stimulating environments for research.  Finally, we thank the referee for their careful reading, comments and suggestions, which helped to clarify the exposition of the paper.

\section{Preliminaries}\label{S:Prelim}
By the usual slight abuse of notation we shall let the set of
morphisms of a small category $\Lambda$ be denoted by $\Lambda$ and identify
an object of $\Lambda$ with its corresponding identity
morphism.  Also note that a monoid $M$ (and hence a group) can be considered as a category
with one object and morphism set equal to $M$, with composition
given by multiplication in the monoid.  For convenience of notation, we shall denote a monoid and its associated category by the same symbol.

The following notation will be used throughout this paper. We let
$\N$ denote the abelian monoid of non-negative integers and we let
$\Z$ be the group of integers.  For a positive integer $k$, we let  $\N^k$
be the product monoid viewed as a category. Similarly, we let $\Z^k$
be the product group viewed, where appropriate, as a category. Let
$\{e_i\}_{i=1}^k$ be the canonical generators of $\N^k$ as a monoid and $\Z^k$ as a group. Moreover, we choose to endow $\N^k$ and $\Z^k$ with the
coordinatewise order induced by the usual order on $\N$ and $\Z$,
i.e. for all $m,n\in\Z^k \; m\le n \iff m-n\in\N^k$.

We will denote by $\K(\mathcal{H})$ the $C^*$-algebra of compact operators on a Hilbert space $\mathcal{H}$.  Where the Hilbert space $\mathcal{H}$ is separable and of infinite dimension we write $\K$ for $\K(\mathcal{H})$.

The concept of a {\it higher rank graph} or $k$-graph
($k=1,2,\ldots$ being the rank) was introduced by A. Kumjian and
D. Pask in \cite{KP00}. We recall their definition of a $k$-graph.
\begin{defn}[{\cite[Definitions 1.1]{KP00}}]
A {\bf $k$-graph} (rank $k$ graph or higher rank graph)
$(\Lambda,d)$ consists of a countable small category $\Lambda$
(with range and source maps $r$ and $s$ respectively) together
with a functor $d:\Lambda\longrightarrow \N^k$ satisfying the
{\bf factorisation property:} for every $\lambda\in\Lambda$ and
$m,n\in\N^k$ with $d(\lambda)=m+n$, there are unique elements
$\mu,\nu\in\Lambda$ such that $\lambda=\mu\nu$ and
$d(\mu)=m,\,d(\nu)=n$. For $n\in\N^k$ and $v\in\Lambda^0$ we
write $\Lambda^n:=d^{-1}(n),\;\Lambda(v):=r^{-1}(v)$ and
$\Lambda^n(v):=\{\lambda\in\Lambda^n \,|\, r(\lambda)=v\}$.
\end{defn}

\begin{defn}[{\cite[Definitions 1.4]{KP00}}]
A $k$-graph $\Lambda$ is {\bf row-finite} if for each
$n\in\N^k$ and $v\in\Lambda^0$ the set $\Lambda^n(v)$ is finite.
We say that $\Lambda$ has {\bf no sources} if
$\Lambda^n(v)\ne\emptyset$ for all $v\in\Lambda^0$ and
$n\in\N^k$.
\end{defn}

Unless stated otherwise, we will assume that each higher rank graph in this paper is row-finite with no sources.  Furthermore, we shall denote such a generic higher rank graph by $(\Lambda,d)$ (or more succinctly $\Lambda$ with the understanding that the degree functor will be denoted by $d$).

We refer to \cite{MacL98} as an appropriate reference on
category theory. There is no need for a detailed knowledge of category
theory as we will be interested in the
combinatorial graph-like nature of higher rank graphs. As the
name suggests a higher rank graph can be thought of as a higher
rank analogue of a directed graph. Indeed, every 1-graph is isomorphic
(in the natural sense)
to the category of finite paths of a directed graph (\cite[Example
1.3]{KP00}). By
\cite[Remarks 1.2]{KP00} $\Lambda^0$ is the set of identity morphisms of $\Lambda$. Indeed it is
fruitful to view $\Lambda^0$ as a set of vertices and $\Lambda$ as a
set of (coloured) paths with composition in $\Lambda$ being concatention of
paths. This viewpoint is discussed further in \cite{E03,RSY03}.

In the sequel we will use the following higher rank graph constructions devised by Kumjian and Pask.  For further examples of $k$-graphs see, for example, \cite{KP00,RSY03,RSY04,S06b}.

\begin{exs}${}$\\[-10pt]\label{exs:k-graph_constructions}
\begin{enumerate}
\item  Let $\Delta_k$ be the category with morphism set equal to $\{(m,n)\in\Z^k\times\Z^k \;|\; m\le n\}$, object set equal to $\{(m,m)\;|\; m\in\Z^k\}$, structure maps defined by $r(m,n)=m,\; s(m,n)=n$ and composition defined by $(m,l)(l,n)=(m,n)$ for all $m,l,n\in\Z^k$.  One may define a degree functor $d:\Delta_k\into\N^k$ by $d(m,n)=n-m$ so that $(\Delta_k, d)$ is a $k$-graph.  Furthermore, it is straightforward to check that $(\Delta_k,d)$ is row-finite and has no sources.

\item \textbf{The product higher rank graph ({\cite[Proposition 1.8]{KP00}})}: Let $(\Lambda_1,d_1)$ and $(\Lambda_2,d_2)$ be rank $k_1,k_2$ graphs respectively, then their product higher rank graph $(\Lambda_1\times\Lambda_2,d_1\times d_2)$ is a $(k_1+k_2)$-graph, where $\Lambda_1 \times \Lambda_2$ is the product category and the degree functor $d_1\times d_2:\Lambda_1\times\Lambda_2 \into \N^{k_1+k_2}$ is given by $d_1\times d_2(\lambda_1,\lambda_2)=(d_1(\lambda_1),d_2(\lambda_2))\in\N^{k_1}\times\N^{k_2}$ for $\lambda_1\in\Lambda_1$ and $\lambda_2\in\Lambda_2$.

\item \textbf{The skew-product higher rank graph ({\cite[Definition 5.1]{KP00}})}: Given a countable group $G$, a $k$-graph $\Lambda$ and a functor
$c:\Lambda\longrightarrow G$, the \emph{skew-product} $k$-graph $G
\times_c \Lambda$ consists of a category with object set identified with $G \times\Lambda^0$ and morphism set identified with $G \times \Lambda$.  The structure maps are given by:
$s(g,\lambda) = (gc(\lambda), s(\lambda))$ and $r(g,\lambda) = (g, r(\lambda))$.
If $s(\lambda)=r(\mu)$ then $(g,\lambda)$ and $(gc(\lambda),\mu)$ are
composable in $G \times_c \Lambda$ and
$(g,\lambda)(gc(\lambda),\mu) = (g,\lambda\mu)$.
The degree map is given by $d(g,\lambda)=d(\lambda)$.  Furthermore, $G$ acts freely on $G\times_c\Lambda$ by $g\cdot(h,\lambda)\mapsto (gh,\lambda)$ for all $g,h\in G$ and $\lambda\in\Lambda$ (see \cite[Remark 5.6]{KP00} and its preceding paragraph).
\end{enumerate}
\end{exs}

To each row-finite $k$-graph with no sources, Kumjian and Pask associated an unique $C^*$-algebra in the following way.

\begin{defn}(\cite[Definitions 1.5]{KP00})
Let $\Lambda$ be a row-finite $k$-graph with no sources.  Then $C^*(\Lambda)$ is defined to be the universal $C^*$-algebra generated by a family $\{s_\lambda \;|\; \lambda\in\Lambda\}$ of partial isometries satisfying:
\begin{enumerate}\renewcommand{\labelenumi}{(\roman{enumi})}
\item $\{s_v \;|\; v\in\Lambda^0\}$ is a family of mutually orthogonal projections,
\item $s_{\lambda\mu}=s_\lambda s_\mu$ for all $\lambda,\mu\in\Lambda$ such that $s(\lambda)=r(\mu)$,
\item $s_\lambda^*s_\lambda = s_{s(\lambda)}$ for all $\lambda\in\Lambda$,
\item for all $v\in\Lambda^0$ and $n\in\N^k$ we have $s_v=\sum_{\lambda\in\Lambda^n(v)} s_\lambda s_\lambda^*$.
\end{enumerate}
For $\lambda\in\Lambda$, define $p_\lambda:=s_\lambda s_\lambda^*$.  A family of partial isometries satisfying (i)--(iv) above is a called a \textbf{$*$-representation} of $\Lambda$.
\end{defn}
We consider the following $C^*$-algebras associated with the constructions noted in Examples \ref{exs:k-graph_constructions}, which will be useful in the sequel.

\begin{exs}${}$\\[-10pt]
\label{exs:k-graph_C*-constructions}
\begin{enumerate}
\item Let $\Delta_k$ be the row-finite $k$-graph with no sources defined in Examples \ref{exs:k-graph_constructions}.1. Then $C^*(\Delta_k)\cong
\K(\ell^2(\Z^k))$ since $\{e_{m,n}\;|\; m,n\in\Z^k\}$ is a complete system of matrix units if $e_{m,n}:=s_{(m,q)}s_{(n,q)}^*$ where $q:=\mbox{sup}\{m,n\}$ (cf. \cite[Examples 1.7 (ii)]{KP00}).
\item Let $(\Lambda_i,d_i)$ be a row-finite $k_i$-graph with no sources for $i=1,2$.  Then $$C^*(\Lambda_1\times \Lambda_2)\cong C^*(\Lambda_1)\otimes C^*(\Lambda_2)$$ by \cite[Corollary 3.5
(iv)]{KP00}.\footnote{The $C^*$-algebra of a row-finite higher rank graph with no sources is nuclear \cite[Theorem 5.5]{KP00}.}
\item Let $G$ be a countable group, $\Lambda$ a row-finite $k$-graph with no sources and
$c:\Lambda\longrightarrow G$ a functor.  Then the action of $G$ on $G\times_c\Lambda$ described in Examples \ref{exs:k-graph_constructions}.3 induces an action $\beta:G\into \aut(C^*(G\times_c\Lambda))$ such that $\beta_g(s_{(h,\lambda)})=s_{(gh,\lambda)}$.  Furthermore $C^*(G\times_c\Lambda)\rtimes_\beta G\cong C^*(\Lambda)\otimes \K(\ell^2(G))$ \cite[Theorem 5.7]{KP00}.
\end{enumerate}
\end{exs}

\section{The $K$-groups of $k$-graph $C^*$-algebras}\label{S:K-theory}
For the remainder of this paper we shall denote by $B_\Lambda$ (or simply $B$ when there is no ambiguity) the $C^*$-algebra of the skew-product of a row-finite $k$-graph $(\Lambda,d)$, with no sources, by $\Z^k$ via the degree functor regarded as a functor into $\Z^k$, i.e. $B:=C^*(\Z^k\times_d\Lambda)$, and by $\beta$ the action of $\Z^k$ on $B$ as described in Examples \ref{exs:k-graph_C*-constructions}.3.  Note that by \cite[Corollary 5.3 and Theorem 5.5]{KP00} and Takesaki-Takai duality \cite{T75} (or \cite[Theorem 5.7]{KP00}, cf. Examples \ref{exs:k-graph_C*-constructions}.3), $C^*(\Lambda)$ is stably isomorphic to the crossed product of an AF-algebra, $B$, by $\Z^k$, i.e. $B\rtimes_\beta \Z^k\cong C^*(\Lambda)\otimes \K$.  Therefore $K_0(C^*(\Lambda))\cong K(B\rtimes_\beta\Z^k)$.

It will be useful for us in the sequel to have an explicit description of how an isomorphism $K_0(C^*(\Lambda))\into K(B\rtimes_\beta\Z^k)$ acts on the $K_0$-class of a canonical projection $p_v,\;v\in\Lambda^0$ in $C^*(\Lambda)$.  To this end, we prefer to investigate how the isomorphism acts by using an alternative approach to that outlined above. This we do below by using a standard technique in $k$-graph $C^*$-algebra theory, namely by utilising the gauge-invariant uniqueness theorem for $k$-graph $C^*$-algebras \cite[Theorem 3.4]{KP00}.

\begin{thrm}\label{thrm:K(C*(Lambda))-K(Crossed-Product)}
Let $\Lambda$ be a row-finite $k$-graph with no sources. Then there exists a group isomorphism $\Phi_0:K_0(B\rtimes_\beta \Z^k)\into K_0(C^*(\Lambda))$ such that $\Phi_0([i_B(p_{(0,v)})]))=[p_v]$ for all $v\in\Lambda^0$ (where we adopt the notation used in \cite{R88} for crossed-product $C^*$-algebras).
\end{thrm}

\begin{proof} \hspace{2pt}Let $(B\rtimes_\beta\Z^k,\,i_B,\,i_{\Z^k})$ be a crossed product for
the dynamical system $(B,\,\Z^k,\,\beta)$ in the sense of \cite{R88}. One checks that $\{t_{(\lambda,(m,n))} \;|\;
  (\lambda,(m,n))\in\Lambda\times\Delta \}$ is a *-representation of
  $\Lambda\times\Delta$, where for
  $(\lambda,(m,n))\in\Lambda\times\Delta$ we let
  $t_{(\lambda,(m,n))}:=i_B(s_{(m,\lambda)})i_{\Z^k}(m+d(\lambda)-n)$.
Moreover $C^*(t_\xi \;|\;
\xi\in\Lambda\times\Delta)=B\rtimes_\beta\Z^k$. Thus by the universal
property of $C^*(\Lambda\times\Delta)$, there exists a $*$-homomorphism
$\pi:C^*(\Lambda\times\Delta)\into B\rtimes_\beta\Z^k$
such that $\pi(s_\xi)=t_\xi$ for all $\xi\in\Lambda\times\Delta$. Let
$\alpha:\T^k\into\aut (B)$ denote the canonical gauge
action on $B$ and let $\hat{\beta}:\T^k\into\aut(B\rtimes_\beta\Z^k)$
denote the dual action of $\beta$. There exists an action
$\tilde{\alpha}$ of $\T^k$
on $B\rtimes_\beta\Z^k$ such that $i_B\alpha_z=\tilde{\alpha}_zi_B$ for
all $z\in\T^k$. It is clear that  setting
$\gamma_{(z_1,z_2)}:=\tilde{\alpha}_{z_1z_2}\hat{\beta}_{z_2^{-1}}$
for all $(z_1,z_2)\in\T^{k}\times\T^k$ defines an action $\gamma$ of
$\T^{2k}$ on $B\rtimes_\beta\Z^k$. Moreover, it satisfies
$\pi\alpha^\times_z=\gamma_z\pi$ for all $z\in\T^{2k}$ where
$\alpha^\times$ is the canonical gauge action on
$C^*(\Lambda\times\Delta)$. Clearly $\pi(p_v)=0$ for all
$v\in\Lambda\times\Delta$, hence by the gauge-invariant uniqueness
theorem \cite[Theorem 3.4]{KP00} we see that
$\pi:C^*(\Lambda\times\Delta)\into B\rtimes_\beta\Z^k$ is a $*$-isomorphism.

For the zero element $0\in\Z^k$, we see that
$p_{(0,0)}$ is a minimal projection in $C^*(\Delta_k)\cong\K$ (cf. Examples \ref{exs:k-graph_C*-constructions}.1).
Therefore, the homomorphism given by $x\mapsto x\otimes p_{(0,0)}$ induces
an isomorphism between $K_0(C^*(\Lambda))$ and
$K_0(C^*(\Lambda)\otimes C^*(\Delta))$, which in turn is isomorphic to
$K_0(C^*(\Lambda\times\Delta)$ (see Examples \ref{exs:k-graph_C*-constructions}.2) and $K_0(B\rtimes_\beta \Z^k)$.  Let $\Psi:K_0(C^*(\Lambda))\into K_0(B\rtimes_\beta\Z^k)$ be the composition of the preceding group isomorphisms.  Then it follows easily that $\Psi([p_v])=[i_B(p_{(0,v)})]$.  Setting $\Phi_0:=\Psi^{-1}$ completes the proof.
\end{proof}

Therefore we may apply \cite[6.10 Theorem]{K88} to describe the $K$-groups
of $C^*(\Lambda)$ by means of a homological spectral sequence with
initial term given by $H_p(\Z^k,K_q(B))$, i.e. the homology of the group $\Z^k$ with coefficients in the left $\Z^k$-module $K_q(B)$ \cite{Br94}, where the $\Z^k$-action is given by $e_i\cdot m = K_0(\beta_{e_i})(m)$ for $i=1,\ldots,k$, (cf. the proof of \cite[Proposition 4.1]{RS01}). First, we recall the definition of a homology spectral sequence from \cite[\S5]{W94} and the notion of convergence (see also \cite{MacL63}).\footnote{The reader will notice that the definition is presented in a less general form than in \cite[\S5]{W94}, but is adequate for our purposes.}

\begin{defn}\label{D:spectral-sequence}
A \emph{homology spectral sequence} (starting at $E^a$) consists of the following data:
\begin{enumerate}
\item A family $\{E^r_{p,q}\}$ of modules defined for all integers $p,q$ and $r\ge a$.
\item Maps $d^r_{pq} : E^r_{p,q} \into E^r_{p-r,q+r-1}$ that are differentials in the sense that $d^r_{p-r,q+r-1} d_{pq}=0$.
\item Isomorphisms $E^{r+1}_{pq}\into \ker(d^r_{pq})/\im (d^r_{p+r,q-r+1})$.
\end{enumerate}
We will denote the above data by $\{(E^r,d^r)\}$.  The \emph{total degree} of the term $E^r_{pq}$ is $n=p+q$.  The homology spectral sequence is said to be \emph{bounded} if for each $n$ there are only finitely many nonzero terms of total degree $n$ in $\{E^r_{pq}\}$, in which case, for each $p$ and $q$ there is an $r_0$ such that $E^r_{pq}\cong E^{r+1}_{pq}$ for all $r\ge r_0$.  We write $E^\infty_{pq}$ for this stable value of $E^r_{pq}$.

We say that a bounded spectral sequence \emph{converges to $\KK_*$} if we are given a family of modules $\{\KK_n\}$, each having a finite filtration
$$
0=F_s(\KK_n)\subseteq \cdots \subseteq F_{p-1}(\KK_n) \subseteq F_p(\KK_n) \subseteq F_{p+1}(\KK_n) \subseteq \cdots \subseteq F_t(\KK_n)=\KK_n,$$
and we are given isomorphisms $E^\infty_{pq} \into F_p(\KK_{p+q})/F_{p-1}(\KK_{p+q})$.
\end{defn}

\begin{lemma}[{cf. \cite[Proposition 4.1]{RS01}}]\label{L:E^infty} There exists a
  spectral sequence $\{(E^r,d^r)\}$ converging to
  $K_*(C^*(\Lambda)):=\{\KK_n\}_{n\in\Z}$ where
$$\KK_n:=\left\{\begin{array}{ll} K_0(C^*(\Lambda)) & \mbox{if $n$ is even},
    \\ K_1(C^*(\Lambda)) & \mbox{if $n$ is odd}.\end{array}\right.$$
Moreover, for $p,q\in\Z$,
    $$E^2_{p,q}\cong \left\{\begin{array}{ll} H_p(\Z^k,K_0(B)) &
        \mbox{if } p\in\{0,1,\ldots,k\} \mbox{ and $q$ is even}, \\ 0
        & \mbox{otherwise},\end{array}\right.$$ $E^\infty_{p,q}\cong
    E^{k+1}_{p,q}$ and $E^{k+1}_{p,q}=0$ if $p\in\Z\backslash\{0,1,\ldots,k\}$ or $q$ is odd.
\end{lemma}

\begin{proof} The first assertion follows from \cite[6.10 Theorem]{K88} applied to
$B\rtimes_\beta\Z^k$, which is $*$-isomorphic to $C^*(\Lambda)\otimes\K$ by Theorem \ref{thrm:K(C*(Lambda))-K(Crossed-Product)}, after noting that $K_*(B \rtimes_\beta \Z^k)$ coincides with its ``$\gamma$-part''
since the Baum-Connes Conjecture with coefficients in an arbitrary
$C^*$-algebra is true for the amenable group $\Z^k$ for all $k\ge 1$.

By the proof of \cite[6.10 Theorem]{K88}, $K_*(B\rtimes_\beta\Z^k)\cong K_*(D)$ for some $C^*$-algebra $D$ which
has a finite filtration by ideals: $0\subset D_0\subset D_1\subset \cdots
\subset D_k=D$ since the dimension of the universal covering space of
the classifying space of $\Z^k$ is $k$.

The spectral sequence we are considering is the spectral sequence $\{(E^r,d^r)\}$ in homology $K_*$
associated with the finite filtration $0\subset D_0\subset D_1\subset
\cdots \subset D_k=D$ of $D$ (\cite[\S 6]{S81}) which has
$E^1_{p,q}=K_{(p+q \mod 2)}(D_p/D_{p-1})$ where $D_n=0$ for $n<0$ and
$D_n=D$ for $n\ge k$. It follows easily that $E^r_{p,q}=0$ for
$p\in\Z\backslash\{0,1,\ldots,k\}$, for all $q\in\Z$ and for all $r\ge
  1$ and $E^\infty_{p,q} \cong E^{k+1}_{p,q}$ (see also \cite[Theorem 2.1]{S81}). This
combined with Kasparov's calculation in the proof of
\cite[6.10 Theorem]{K88}, giving $E^2_{p,q}\cong H_p(\Z^k,K_q(B))$,
along with the observation that $K_q(B)=0$ for odd $q$, since $B$ is an
AF-algebra, proves the second assertion.\end{proof}

Now we will compute $H_*(\Z^k,K_0(B))$ in terms of the combinatorial data encoded in $\Lambda$.
First, let us examine the structure of $B$, and hence $K_0(B)$, in a little more detail.

\begin{lemma}\label{L:structure_of_B}
Let $\Lambda$ be a row-finite $k$-graph with no sources.  Then
$$
B =  \overline{\bigcup_{n\in\Z^k} B_n},
$$where
 $$B_n\!=\!\Span\{s_\lambda
    s_\mu^* | \lambda,\mu\in \Z^k \times_d \Lambda,
    s(\lambda)\!=\!s(\mu)\!=\!(n,v)\;\mbox{for some}\;v\in\Lambda^0 \}\cong\bigoplus_{v\in\Lambda^0} \!B_n(v),$$ and
$$B_n(v):=\Span\{s_\lambda
    s_\mu^* \;|\; \lambda,\mu\in \Z^k \times_d \Lambda,\;
    s(\lambda)=s(\mu)=(n,v) \}\cong \K(\ell^2(s^{-1}(v)))$$ for all $v\in
    \Lambda^0$ and $n\in\Z^k$.
\end{lemma}

\begin{proof}
Follows immediately from the proofs of \cite[Lemma 5.4, Theorem 5.5]{KP00} and the observation that for all $n\in\Z^k$ and $v\in\Lambda^0$, $s^{-1}((n,v))\subset \Z^k\times_d\Lambda$ may be identified with $s^{-1}(v)\subset\Lambda$ via $(n-d(\lambda),\lambda)\mapsto
    \lambda$ for all $\lambda\in s^{-1}(v)$.
\end{proof}

\begin{defn}\label{D:ZLambda^0}
Let $\Z\Lambda^0$ be the group of all maps from $\Lambda^0$ into $\Z$ that have finite support under pointwise addition.  For each $u\in\Lambda^0$, we denote by $\delta_u$ the element of $\Z\Lambda^0$ defined by $\delta_u(v)=\delta_{u,v}$ (the Kronecker delta) for all $v\in\Lambda^0$.  Note that $\Z\Lambda^0$ is a free abelian group with free set of generators $\{\delta_u\;|\; u\in\Lambda^0\}$.
\end{defn}
\begin{defn}[cf. {\cite[\S6]{KP00}}]
      Define the vertex matrices of $\Lambda$, $M_i$, by the
      following. For $u,v\in\Lambda^0$ and $i=1,2,\ldots,k$,
      $M_i(u,v):=|\{\lambda\in\Lambda^{e_i}\,|\,r(\lambda)=u,
      s(\lambda)=v\}|$.
    \end{defn}
\begin{rem}\label{R:vertex_mxs_commute}
By the factorisation property, the vertex matrices of a $k$-graph pairwise commute \cite[\S6]{KP00}.
\end{rem}

\begin{convention}\label{C:matrix-endo}Given a matrix $M$ with integer entries and index set $I$, by slight abuse of notation we shall, on occasion, regard $M$ as the group endomorphism $\Z\Lambda^0\into \Z\Lambda^0$, defined in the natural way as $(Mf)(i)=\sum_{j\in I}M(i,j)f(j)$ for all $i\in I,\;f\in\Z\Lambda^0$.
\end{convention}

    \begin{lemma}[cf. \cite{RS01} Lemma 4.5 and \cite{PR96} Proposition
      4.1.2]\label{L:K_0(B)} For all $n,m\in\Z^k$ such that $m\le n$, let $A_m:=\Z\Lambda^0$.  Moreover, define homomorphisms
  $\jmath_{nm}:A_m\into A_n$ by $\jmath_{mm}(f)=f$,
  $\jmath_{m+e_i, m}:=M_i^t$ for all $f\in A_m,\; u\in\Lambda^0,\;i\in\{1,\ldots,k\}$, and $\jmath_{m+e_i+e_j,m}:=\jmath_{m+e_i+e_j,m+e_i}\jmath_{m+e_i,m}$ for all $j\in\{1,\ldots,k\}$. Then $(A_m;\jmath_{nm})$ is a direct system of groups and $K_0(B)$ and
  $A:=\lim_{\rightarrow}(A_{m};\jmath_{nm})$ are isomorphic.
    \end{lemma}

    \begin{proof} It follows from Remarks \ref{R:vertex_mxs_commute} that the connecting homomorphisms are well-defined and that $(A_m;\jmath_{nm})$ is a direct system. From Lemma \ref{L:structure_of_B} we deduce that $\displaystyle K_0(B)
    \cong \lim_{\rightarrow} (K_0(B_n); K_0(\iota_{n,m}))$, where, for
    $m,n\in\Z^k$ with $m\le n$, $\iota_{n,m}: B_m \into B_n$ are the
    inclusion maps \cite[Proposition 6.2.9]{W-O93}.  We also deduce that for $n\in\Z^k$ and $v\in
    \Lambda^0$, $K_0(B_n(v))$ is isomorphic to $\Z$ and is generated
    by the equivalence class comprising all minimal projections in $B_n(v)$, of which $p_\xi$ is a member for any $\xi\in s^{-1}(n,v)$. Therefore, $K_0(B_n)
    \cong \bigoplus_{v\in \Lambda^0}K_0(B_n(v))$ is generated by $\{[p_{(n,v)}]_n \;|\; v\in\Lambda^0\}$, where
    $[\,\cdot\,]_n$ denotes the equivalence classes of $K_0(B_n)$ for all
    $n\in \Z^k$. Thus the map $\psi_n: A_n
    \longrightarrow K_0(B_n)$ given by $f\mapsto \sum_{u\in\Lambda^0} f(u) [p_{(n,u)}]_n$
    is a group isomorphism for all $n\in\Z^k$.

The embedding
    $\iota_{n,m}:B_m \longrightarrow B_n$ is given by
    \[
    \iota_{n,m}(s_{(m-d(\lambda),\lambda)} s_{(m-d(\mu),\mu)}^*) =
    \sum_{\alpha\in\Lambda^{n-m}(s(\lambda))}
    s_{(m-d(\lambda),\lambda\alpha)} s_{(m-d(\mu),\mu\alpha)}^*
    \]
    for all $\lambda,\mu\in\Lambda$. Therefore,
    \begin{eqnarray*}
      K_0(\iota_{n+e_i, n}) \left( [p_{(n,v)}]_n \right) &=& \left[
        \iota_{n+e_i, n} \left( p_{(n,v)} \right)
      \right]_{n+e_i}
      = \left[ \sum_{\alpha\in\Lambda^{e_i}(v)} p_{(n,\alpha)}
      \right]_{n+e_i}\\
      &=& \sum_{u\in\Lambda^0} M_i(v,u)[p_{(n+e_i,u)}]_{n+e_i}
    \end{eqnarray*}
    and
    \begin{eqnarray*}
      K_0(\iota_{n+e_i, n})\left(\sum_{v\in\Lambda^0} f(v)
        [p_{(n,v)}]_n\right) &=& \sum_{u\in\Lambda^0} \left(
        \sum_{v\in\Lambda^0} M_i(v,u)f(v) \right)
      [p_{(n+e_i,u)}]_{n+e_i}\\
      &=& \sum_{u\in\Lambda^0} \jmath_{n+e_i,n}(f)(u)[p_{(n+e_i,u)}]_{n+e_i}.
    \end{eqnarray*}

    Thus the following square commutes for all $i\in\{1,2,\ldots,k\}$.
    \[
    \begin{CD} K_0(B_n) @>{K_0( \iota_{n+e_i, n}
        )}>> K_0(B_{n+e_i})\\
      @A{\psi_n}AA @AA{\psi_{n+e_i}}A\\
      A_n @>{\jmath_{n+e_i,n}}>> A_{n+e_i}
    \end{CD}
    \]

    The result follows.\end{proof}

Henceforth, we shall follow the notation introduced in Lemma \ref{L:K_0(B)} and its proof.

Now we begin to examine the action of $\Z^k$ on $K_0(B)$ in terms of the description of $K_0(B)$ provided by Lemma \ref{L:K_0(B)}.

\begin{lemma}[cf. \cite{RS01} Lemma 4.10]\label{L:K_0(beta)} Fix $i\in\{1,\ldots,k\}$ and
  define a homomorphism $\phi_{i,n}:A_n\into A_n$ by $\phi_{i,n}:=M_i^t$ for all $n\in\Z^k$. Let $\phi_i: A \into A$ be the homomorphism induced by the system of homomorphisms $\{\phi_{i,n}\;|\; n\in\Z^k\}$. Then
  $\psi\phi_i=K_0(\beta_{e_i})\psi$, where $\psi: A \into K_0(B)$ is the isomorphism constructed in
  Lemma \ref{L:K_0(B)}.
\end{lemma}

\begin{proof} It follows from Remarks \ref{R:vertex_mxs_commute} that $\phi_{i,n}\jmath_{nm}=\jmath_{nm}\phi_{i,m}$ for all $m,n\in\Z^k$ so that $\phi_i$ is well-defined for all $i\in\{1,\ldots,k\}$.

Now, we let $\tilde{\psi}: A \into \lim_{\rightarrow} (K_0(B_m);
K_0(\iota_{n,m}))$ be the unique isomorphism such that
$K_0(\iota_n)\circ \psi_n = \tilde{\psi} \circ \jmath_n$ for all $n\in\Z^k$
where $\psi_n:A_n \into K_0(B_n): f \mapsto \sum_{u\in\Lambda^0}
f(u)[p_{(n,u)}]_n$ (cf. proof of Lemma \ref{L:K_0(B)}). Then $\psi:
A\into K_0(B)$ is the composition of $\tilde{\psi}$ with the canonical
isomorphism of $\lim_\rightarrow (K_0(B_n); K_0(\iota_{n,m}))$ onto $K_0(B)$.

We will show that
\[
\begin{CD} K_0(B_n) @>{K_0(\iota_n)}>> K_0(B) \\
@V{\tilde{\phi}_{i,n}}VV  @VV{K_0(\beta_{e_i})}V\\
K_0(B_n) @>{K_0(\iota_n)}>> K_0(B)
\end{CD}
\]
commutes for all $i=1,2,\ldots,k$ and $n\in\Z^k$ where
$\iota_n:B_n \into B$ is the inclusion map for all $n\in\Z^k$ and
$\tilde{\phi}_{i,n}=\psi_n\circ\phi_{i,n}\circ \psi_n^{-1}$. For then the Lemma
follows from the universal properties of direct limits.

Fix $i\in\{1,\ldots,k\}$ and $n\in\Z^k$. Since $K_0(B_n)$ is generated by
$\{[p_{(n,v)}]_n \;|\; v\in\Lambda^0 \}$, it suffices to show
that $K_0(\beta_{e_i})\circ K_0(\iota_n)([p_{(n,v)}]_n) =
K_0(\iota_n)\circ \tilde{\phi}_{i,n} ([p_{(n,v)}]_n)$ for all $v\in\Lambda^0$. To
see that this holds let $v\in\Lambda^0$, then
\[ K_0(\beta_{e_i})\circ K_0(\iota_n)([p_{(n,v)}]_n) =
K_0(\beta_{e_i})([p_{(n,v)}]) = [p_{(n+e_i,v)}]. \]
While \[ K_0(\iota_n)\circ \tilde{\phi}_{i,n} ([p_{(n,v)}]_n) =
\sum_{u\in\Lambda^0}M_i(v,u)[p_{(n,u)}] =
\sum_{\alpha\in\Lambda^{e_i}(v)}[p_{(n+e_i,\alpha)}] = [p_{(n+e_i,v)}].
\]\end{proof}

Having established a description of $K_0(B)$ as a left $\Z^k$-module in terms of the structure of $\Lambda$, we are almost in a position to describe $H_*(\Z^k,K_0(B))$.  First we recall some relevant notions from homological algebra.

It will be convenient to use multiplicative notation for the free abelian group $\Z^k$, generated by $k$ generators. Thus we set $G:= \langle s_i \;|\; s_i s_j=s_js_i \mbox{ for all } i,j\in\{1,\ldots,k\} \rangle$ and $R:=\Z G$, the group ring of $G$ \cite{Br94}. An efficient method of computing $H_*(\Z^k,K_0(B))$ is by means of a
Koszul resolution $K(\mbf{x})$ for an appropriate regular sequence $\mbf{x}$ on $R$, \cite[Corollary 4.5.5]{W94}.

By a regular sequence on $R$ we mean a sequence $\mbf{x}=\{x_i\}_{i=1}^n$ of elements of $R$ such that
\begin{itemize}
\item[(a)] $(x_1,\ldots,x_n)R \ne R$; and
\item[(b)] For $i=1,\ldots,n,\; x_i\not\in
  \mathcal{Z}(R/(x_1,\ldots,x_{i-1})R)$.
\end{itemize}
In statements (a) and (b), we regard $R$ as an $R$-module; denote by $(x_1,\ldots,x_j)\;(j=1,\ldots,n)$ the ideal of $R$ generated by $\{x_i\}_{i=1}^j$ and $(x_1,\ldots,x_j)R$ the sub-$R$-module $\{r\cdot r' \;|\; r\in (x_1,\dots,x_j), r' \mbox{ in the $R$-module } R\}$ of $R$; and denote by $\mathcal{Z}(M)$ the set of zero-divisors on a R-module $M$, i.e. $\mathcal{Z}(M):=\{ r\in R \;|\; r\cdot m =0$ for some
  non-zero $m\in M \}$ (see \cite[\S 3.1]{Kap74} for more details).

It is straightforward to check that for any finite set of generators $\{t_1,\ldots,t_k\}$ of $G$, the subset $\mbf{x}=\{x_i\}_{i=1}^k$, where $x_i:=1-t_i$ for $i=1,\ldots,k$, is a regular sequence on $R$.

Following \cite[\S4.5]{W94} we will describe the Koszul complex, $K(\mbf{x})$, in terms of the exterior algebra of a free $R$-module \cite{AB74}.  It will be convenient for us to describe the terms of the exterior algebra as follows.

\begin{defn} For any non-negative integer $l$, let $\E^l(R^k)$
denote the $l^{\mbox{\scriptsize th}}$ term of the exterior algebra of the free $R$-module
$R^k:=\bigoplus_{i=1}^k R_i$, over
$R$, where $R_i=R$ for $i=1,\ldots,k$. Moreover, for any negative integer $l$, let $\E^l(R^k)=\{0\}$.

{\setlength{\parindent}{0cm}
For $l\in\Z$ let $N_l:= \left\{\!\begin{array}{ll} \{
    (\mu_1,\ldots,\mu_l)\in \{1,\ldots,k\}^l \;|\;
    \mu_1<\cdots<\mu_l \} & \!\!\mbox{if } l\in\{1,\ldots,k\},\\ \{\star\} &
    \!\!\mbox{if } l=0,\\
    \emptyset & \!\!\mbox{otherwise.}\end{array}\right.$

For $l\in\{1,\ldots,k\}$ and $\mu=(\mu_1,\ldots,\mu_l)\in N_l, i=1,\ldots,l$ we let
    $$\mu^i:=\left\{\begin{array}{ll}(\mu_1,\mu_2\ldots,\mu_{i-1},\mu_{i+1},\mu_{i+2}\ldots,\mu_{l})\in
    N_{l-1} & \mbox{if } l\ne 1,\\ \star & \mbox{if } l=1.\end{array}\right.$$

For $n=1,2,\dots$ and $r\in\Z$, let $${n \choose r }:= \left\{\begin{array}{ll}
    \frac{n!}{(n-r)!r!} & \mbox{if } 0 \le r \le n, \\ 0 & \mbox{if }
    r<0 \mbox{ or } r>n.\end{array}\right.$$

Using the above notation we may describe the $l^{\mbox{\scriptsize th}}$-term ($l\in\Z$)  of the exterior algebra over $R^k$, $\E^l(R^k)$, as being generated by the set $N_l$ as a free $R$-module and having rank $k \choose l$.} \end{defn}

Now let $K(\mbf{x})$ be the chain complex
\[ 0 \longleftarrow \E^0(R^k) \longleftarrow \E^1(R^k) \longleftarrow
\cdots \longleftarrow \E^k(R^k) \longleftarrow 0 \]
where the differentials $\E^l(R^k)\into \E^{l-1}(R^k)$ are given by mapping
\[ \mu \mapsto \sum_{j=1}^l (-1)^{j+1}
x_{\mu_j}\mu^j \qquad \mbox{for all } \mu=(\mu_1,\ldots,\mu_l)\in
N_l,\] if $l\in\{1,\ldots,k\}$ and the zero map otherwise. By \cite[Corollary 4.5.5]{W94} $K(\mbf{x})$ is a free
resolution of $R/I$ over $R$ where $I$ is the ideal of $R$
generated by $\mbf{x}$. It is well known (see e.g. \cite[Chapter
6]{W94},\cite[\S I.2]{Br94}) that $I=\ker \epsilon$ where
$\epsilon:R\into \Z: g \mapsto 1$ is the augmentation map of the
group ring $R=\Z G$. Thus we have a free (and hence projective) resolution of $\Z$ over
$\Z G$, which we may use to compute $H_*(G, K_0(B))$ (see \cite[Chapter III]{Br94}).

\begin{lemma}\label{L:Koszul}
Following the above notation, we have $H_*(G,K_0(B))$ isomorphic to the homology of the chain complex
\begin{equation*} \mathcal{B}: 0 \longleftarrow K_0(B) \longleftarrow \cdots
    \longleftarrow \bigoplus_{N_l} K_0(B) \longleftarrow
    \cdots \longleftarrow K_0(B) \longleftarrow 0, \end{equation*} where the
    differentials $\tilde{\partial_l}:\bigoplus_{N_l} K_0(B) \into
    \bigoplus_{N_{l-1}}K_0(B)\; (l\in\{1,\ldots,k\})$ are defined by
    \[ \bigoplus_{\mu\in N_{l}} m_\mu \mapsto \bigoplus_{\lambda\in
      N_{l-1}}
    \sum_{\mu\in N_l} \sum_{i=1}^l (-1)^{i+1}
    \delta_{\lambda,\mu^i}(m_\mu-K_0(\beta_{\mu_i})(m_\mu)). \]
    (Recall that the $G$-action on $K_0(B)$ is given by
    $s_i\cdot m=K_0(\beta_{e_i})(m)$ for all $m\in
    K_0(B),\;i=1,\dots,k$.)
\end{lemma}

\begin{proof}
By definition $H_*(G,K_0(B))\cong H_*(K(\mbf{x})\otimes_G K_0(B))$, where the latter chain complex is obtained by applying the functor $-\otimes_G K_0(B)$ termwise to the chain complex $K(\mbf{x})$.  The Lemma follows from the fact that $\E^l(R^k)\otimes_G K_0(B)$ is canonically isomorphic to $\bigoplus_{N_l} K_0(B)$ ($l\in\{1,\ldots,k\}$) and setting $t_i:=s_i^{-1}$ ($i=1,\ldots,k$) as our generators of $G$ to obtain $\mbf{x}$ as described above.
\end{proof}

For $m,n\in\Z^k$ with $m\le n$, let $\mathcal{A}^{(n)}$ be the chain complex
$$  0 \longleftarrow A_n \longleftarrow \cdots
    \longleftarrow \bigoplus_{N_l} A_n \longleftarrow
    \cdots \longleftarrow A_n \longleftarrow 0, $$ with
    $A_n=\Z\Lambda^0$ and differentials,
$\partial^{(n)}_l:\bigoplus_{N_l} A_n \into
    \bigoplus_{N_{l-1}} A_n \;(l\in\{1,\ldots,k\})$, defined by
    \[ \bigoplus_{\mu\in N_{l}} m_\mu \mapsto \bigoplus_{\lambda\in
      N_{l-1}}
    \sum_{\mu\in N_l} \sum_{i=1}^l (-1)^{i+1}
    \delta_{\lambda,\mu^i}(m_\mu - \phi_{\mu_i,n}(m_\mu)), \]
where for $i=1,\ldots,k$ and $n\in\Z^k$, $\phi_{i,n}$ is the homomorphism defined in Lemma \ref{L:K_0(beta)}.
Furthermore, let $(\tau_m^n)_p: \mathcal{A}^{(m)}_p
\into \mathcal{A}^{(n)}_p$ be the homomorphism
defined by $(\tau_m^n)_p(\bigoplus_{\mu\in N_p} m_\mu) = \bigoplus_{\mu\in
  N_p} \jmath_{nm}(m_\mu)$ for all $p\in\{0,1,\ldots,k\}$ and the trivial map for $p\in\Z\backslash\{0,1,\ldots,k\}$ (cf. Lemma
\ref{L:K_0(B)}).

Following \cite[Chapter 4, \S1]{S66}, by a chain map $\tau:\mathcal{C}\into \mathcal{C}'$ we mean a collection $\{\tau_p:\mathcal{C}_p\into \mathcal{C}'_p\}$ of homomorphisms that commute with the differentials in the sense that commutativity holds in each square:
$$
\begin{CD}
\mathcal{C}_p @>>> \mathcal{C}_{p-1} \\
@V\tau_pVV  @VV\tau_{p-1}V\\
\mathcal{C}'_p @>>> \mathcal{C}'_{p-1}
\end{CD}
$$
Recall that there is a category of chain complexes whose objects are chain complexes and whose morphisms are chain maps.  Moreover, the category of chain complexes admits direct limits.

\begin{lemma}
Following the above notation, for each $m,n\in\Z^k$ the system of homomorphisms $\{(\tau_m^n)_p\;|\; p\in\Z\}$ defines a chain map $\tau_m^n:\mathcal{A}^{(m)}\into\mathcal{A}^{(n)}$ such that $(\mathcal{A}^{(n)}; \tau_m^n)$ is a direct system in the category of chain complexes. Furthermore, $(\mathcal{B};\gamma_n)$ is a direct limit for $(\mathcal{A}^{(n)};\tau_m^n)$, where $\gamma_n:\mathcal{A}^{(n)}\into \mathcal{B}$ is given by $(\gamma_n)_p(\bigoplus_{\mu\in N_p} m_\mu)=\bigoplus_{\mu\in N_p}\psi\jmath_n(m_\mu)$ for all $p\in\{0,1,\ldots,k\}$ and the trivial map otherwise.
\end{lemma}
\begin{proof}
That $\tau_m^n : \mathcal{A}^{(m)}\into\mathcal{A}^{(n)}$ is a chain map for all $m,n\in\Z^k$, follows immediately from the fact that $\phi_{i,n}\jmath_{nm}=\jmath_{nm}\phi_{i,m}$ for all $i\in\{1,\ldots,k\}$ and $m,n\in\Z^k$ (cf. proof of Lemma \ref{L:K_0(beta)}).  That $(\mathcal{A}^{(n)}; \tau_m^n)$ is a direct system of chain complexes follows immediately from the fact that $(A_m;\jmath_{nm})$ is a direct system of groups (by Lemma \ref{L:K_0(B)}).

Note that since $K_0(\beta_{e_i})\psi=\psi\phi_i$ for all $i\in\{1,\dots,k\}$ by Lemma \ref{L:K_0(beta)}, a direct calculation shows that $\tilde{\partial}_p(\gamma_n)_p=(\gamma_n)_{p-1}\partial^{(n)}_p$ for all $p\in\Z^k$, thus $\gamma_n:\mathcal{A}^{(n)}\into \mathcal{B}$ is a chain map for all $n\in\Z^k$.  The fact that $\gamma_m=\gamma_n\tau_{nm}$ for all $m,n\in\Z^k$ follows immediately by construction of the maps.

Now suppose that $(\mathcal{A};\tau_n)$ is a direct limit for $(\mathcal{A}^{(n)};\tau_m^n)$.  Then by the above and the universal property of direct limits, there exists a morphism $\gamma:\mathcal{A} \into \mathcal{B}$ such that $\gamma \tau_m^n = \gamma_n$ for all $m,n\in\Z^k$.  In order to show that $\gamma$ is an isomorphism it suffices to show that each $(\gamma)_p:\mathcal{A}_p\into \mathcal{B}_p$ is an isomorphism for all $p\in\Z$.  This follows immediately from the fact that $\psi:A\into K_0(B)$ is an isomorphism (Lemma \ref{L:K_0(B)}) and that direct limits commute with (finite) direct sums in the category of abelian groups in the obvious way.  We have shown therefore that $(\mathcal{B};\gamma_n)$ is a direct limit for $(\mathcal{A}^{(n)};\tau_m^n)$ in the category of chain complexes.
\end{proof}

Note that each chain complex $\mathcal{A}^{(n)}$ does not actually
  depend on $n\in\Z^k$, thus for ease of notation we let $\mathcal{D}$ denote this common
  chain complex with differentials $\partial_p:=\partial^{(n)}_p$ for
  all $p\in\Z$.

\begin{thrm}\label{T:Hom} Using the above notation, the homology of $\Z^k$ with coefficients in the left $\Z^k$-module $K_0(B)$ is given by the homology of the chain complex $\mathcal{D}$, i.e. we have $H_*(G,K_0(B))\cong \Homo_*(\mathcal{D})$.
\end{thrm}

\begin{proof} The homology functor commutes with direct limits (\cite[Chapter 4,
\S 1, Theorem 7]{S66}), therefore it follows that $H_*(G,K_0(B))\cong
\lim_{\rightarrow}( \Homo_*(\mathcal{A}^{(n)}), \Homo_*(\tau_m^n) )$. Thus, it suffices to prove that $\Homo(\tau_m^{m+e_j})_p$ is the
identity map for all $p\in\Z,\;m\in\Z^k,\;j\in\{1,\ldots,k\}$. To see that this
is true we show that
$\bigoplus_{\mu\in N_p}(1-M_j^t)(y)\in \im \partial_{p+1}$ for all
$y\in \ker\partial_{p},\;p\in\Z,\;j\in\{1,\ldots,k\}$. Indeed, we claim that given
$y=\bigoplus_{\mu\in N_p} y_\mu \in \ker \partial_p$ we have
$$ \bigoplus_{\mu\in N_p} (1-M_j^t)y_\mu =
\partial_{p+1}\left(\bigoplus_{\lambda\in N_{p+1}}
  z_{\lambda}\right)$$ where $z_{\lambda}=\sum_{i=1}^{p+1} (-1)^{i+1}
\delta_{\lambda_i,j} y_{\lambda^i}$ for all $\lambda\in N_{p+1}$.

Fix $j,p\in\{1,\ldots,k\}$. Let $\kappa(\mu):=\sum_{i=1}^p\delta_{j,\mu_i}i$, i.e. if
$j$ is a component of $\mu$, then $\kappa(\mu)$ denotes the unique
$i\in\{1,\ldots,k\}$ such that $\mu_i=j$; otherwise $\kappa(\mu)=0$. Now
fix a $\mu'=(\mu_1',\dots,\mu_p')\in N_p$ and let
$y=\bigoplus_{\mu\in N_p} y_\mu$ be in $\ker \partial_p$.

First, suppose that $\kappa(\mu')>0$ and let $\eta=(\mu')^{\kappa(\mu')}$.
Then
\begin{eqnarray*} 0 &=& \partial_p(y)_\eta = \sum_{\mu\in N_p}
  \sum_{i=1}^p (-1)^{i+1} \delta_{\eta,\mu^i} (1-M_{\mu_i}^t)
  y_{\mu} \\ &=& \sum_{i=1}^p (-1)^{i+1}
  \delta_{\eta,(\mu')^i}(1-M_{\mu_i}^t)x_{\mu'} + \sum_{\mu\in N_p
    \atop \mu\ne\mu'} \sum_{i=1}^p (-1)^{i+1}
  \delta_{\eta,\mu^i}(1-M_{\mu_i}^t)x_\mu \\ &=&
  (-1)^{\kappa(\mu')+1}(1-M_j^t)x_{\mu'}  + \sum_{\mu\in N_p
    \atop \mu\ne\mu'} \sum_{i=1}^p (-1)^{i+1}
  \delta_{\eta,\mu^i}(1-M_{\mu_i}^t)x_\mu,
  \end{eqnarray*}so that
\begin{eqnarray*}
  (1-M_j^t)x_{\mu'} &=& \sum_{\mu\in N_p \atop \mu\ne \mu'} \sum_{i=1}^p(-1)^{i+\kappa(\mu')+1}\delta_{\eta,\mu^i}(1-M_{\mu_i}^t)x_\mu.\end{eqnarray*}
Now \begin{eqnarray*} \partial_{p+1}\left(\bigoplus_{\lambda\in
      N_{p+1}} z_{\lambda}\right)_{\mu'} &=& \sum_{\lambda\in N_{p+1}}
    \sum_{i,r=1}^{p+1} (-1)^{i+r+2}
    \delta_{\mu',\lambda^i}\delta_{j,\lambda_r}(1-M_{\lambda_i}^t)y_{\lambda^r}\\
    &=& \sum_{\lambda\in N_{p+1} \atop \kappa(\lambda)>0} \sum_{i=1}^{p+1}
    (-1)^{i+\kappa(\lambda)}
    \delta_{\mu',\lambda^i}(1-M_{\lambda_i}^t)y_{\lambda^{\kappa(\lambda)}}
    \\
&=& \sum_{\lambda\in N_{p+1} \atop \kappa(\lambda)>0}
\left\{\sum_{i=1}^{\kappa(\lambda)-1} (-1)^{i+\kappa(\mu')+1}
\delta_{\mu',\lambda^i}(1-M_{\lambda_i}^t)y_{\lambda^{\kappa(\lambda)}}\right. \\ &+&
\left.\sum_{i=\kappa(\lambda)+1}^{p+1} (-1)^{i+\kappa(\mu')}
  \delta_{\mu',\lambda^i}(1-M_{\lambda_i}^t)y_{\lambda^{\kappa(\lambda)}} \right\}\\
\end{eqnarray*}
\begin{eqnarray*}
&=& \sum_{\mu\in N_p \atop \kappa(\mu)=0} \sum_{i=1}^p
(-1)^{i+\kappa(\mu')+1} \delta_{\eta,\mu^i}(1-M_{\mu_i}^t)y_\mu
\end{eqnarray*}

since for every $\lambda\in N_{p+1}$ such that $\kappa(\lambda)>0$ and for
every $i\in\{1,\dots,p+1\}\backslash\{\kappa(\lambda)\}$ we have \begin{enumerate} \item $\delta_{\mu',\lambda^{\kappa(\lambda)}}=0$,
\item
  $\kappa(\lambda)=\left\{\begin{array}{ll} \kappa(\mu')+1 & \mbox{if }
      \mu'=\lambda^i \mbox{ with } i < \kappa(\lambda), \\
      \kappa(\mu') & \mbox{if } \mu'=\lambda^i \mbox{ with } \kappa(\lambda) < i,\end{array}\right.$
\item $\mu' = \lambda^i \iff \eta = \left\{\begin{array}{ll}
      (\lambda^{\kappa(\lambda)})^i & \mbox{if } i<\kappa(\lambda), \\
      (\lambda^{\kappa(\lambda)})^{i-1} & \mbox{if }
      \kappa(\lambda)<i,\end{array}\right.$
\end{enumerate}
and if $\mu\in N_p$ then $\kappa(\mu)=0,\;\eta=\mu^i$ for some $i\in\{1,\ldots,k\} \iff
\mu\ne \mu',\; \eta=\mu^i$ for some $i\in\{1,\dots,k\}$.

Hence, \begin{eqnarray*} \partial_{p+1}\left( \bigoplus_{\lambda\in N_{p+1}}
  z_\lambda \right)_{\mu'} &=& \sum_{\mu\in N_p \atop \mu\ne \mu'}
\sum_{i=1}^p (-1)^{i+\kappa(\mu')+1}\delta_{\eta,\mu^i}(1-M_{\mu_i}^t)y_\mu \\
  &=& (1-M_j^t) y_{\mu'}.\end{eqnarray*}

Now suppose that $\kappa(\mu')=0$. Then
\begin{eqnarray*} \partial_{p+1}\left( \bigoplus_{\lambda\in N_{p+1}}
    z_\lambda \right)_{\mu'} &=& \sum_{\lambda\in N_{p+1}}
  \sum_{i,r=1}^{p+1}
  (-1)^{i+r+2}\delta_{\mu',\lambda^i}\delta_{j,\lambda_r}(1-M_{\lambda_i}^t)y_{\lambda^r}
\\ &=& (-1)^{\kappa(\xi)+\kappa(\xi)+2}(1-M_{\xi_{\kappa(\xi)}}^t)y_{\xi^{\kappa(\xi)}}
\\ &=& (1-M_j^t)y_{\mu'},
\end{eqnarray*}
where $\xi$ is the unique element of $N_{p+1}$ satisfying $\kappa(\xi)>0$
and $\xi^{\kappa(\xi)}=\mu'$.\end{proof}

Combining the results of this section we get the following theorem.

\begin{thrm}\label{T:Main} Let $\Lambda$ be a row-finite $k$-graph with no
  sources. Then there exists a spectral sequence $\{(E^r,d^r)\}$ converging to
  $K_*(C^*(\Lambda))$ with $E^\infty_{p,q}\cong E^{k+1}_{p,q}$ and
  $$E^2_{p,q}\cong \left\{\begin{array}{ll}\Homo_p(\mathcal{D}) &
      \mbox{if } p\in\{0,1,\ldots,k\} \mbox{ and $q$ is even,}\\ 0 & \mbox{otherwise,}\end{array}\right.$$ where
  $\mathcal{D}$ is the chain complex with
  $$\mathcal{D}_p:=\left\{\begin{array}{ll} \bigoplus_{\mu\in N_p}\Z\Lambda^0 & \mbox{if }p\in\{0,1,\ldots,k\},\\
  0 & \mbox{otherwise.}\end{array}\right. $$
  and differentials $$\partial_p:\mathcal{D}_p\into\mathcal{D}_{p-1}:
 \bigoplus_{\mu\in N_{p}} m_\mu \mapsto \bigoplus_{\lambda\in
      N_{p-1}}
    \sum_{\mu\in N_p} \sum_{i=1}^p (-1)^{i+1}
    \delta_{\lambda,\mu^i}(1-M_{\mu_i}^t)m_\mu$$ for $p\in\{1,\ldots,k\}$.
\end{thrm}

Specialising Theorem \ref{T:Main} to the case when $k=2$ gives us
explicit formulae to compute the $K$-groups of the $C^*$-algebras of row-finite $2$-graphs with no sources.

\begin{prop}\label{P:k=2} Let $\Lambda$ be a row-finite 2-graph with no sources and vertex matrices $M_1$ and $M_2$.  Then there is an isomorphism
$$
\Phi:\coker(1-M_1^t,1-M_2^t)\oplus\ker{ M_2^t-1 \choose 1-M_1^t }\into K_0(C^*(\Lambda))
$$
such that $\Phi((\delta_u + \im\partial_1)\oplus 0)=[p_u]$ for all $u\in\Lambda^0$ (cf. Definition \ref{D:ZLambda^0}) and where we regard $(1-M_1^t,1-M_2^t):\Z\Lambda^0\oplus\Z\Lambda^0\into
\Z\Lambda^0$ and $\displaystyle{ M_2^t-1 \choose 1-M_1^t }:\Z\Lambda^0 \into
\Z\Lambda^0\oplus\Z\Lambda^0$ as group homomorphisms defined in the
natural way.

Moreover, we have
$$
K_1(C^*(\Lambda)) \cong
\ker(1-M_1^t,1-M_2^t)/\im{ M_2^t-1 \choose 1-M_1^t }.
$$
\end{prop}

\begin{proof} The Kasparov spectral sequence converging to $K_*(C^*(\Lambda))$ of Proposition
\ref{T:Main} has $E^\infty_{p,q}\cong E^3_{p,q}$ for all
$p,q\in\Z$. However, it follows from $E^2_{p,q}=0$ for
odd $q$ that the differential $d^2$ is
the zero map and $E^3_{p,q}\cong E^2_{p,q}\cong \Homo_p(\mathcal{D})$ for
all $p\in\{0,1,\ldots,k\}$ and even $q$, where $\mathcal{D}$ is the chain
complex
$$0 \outof \Z\Lambda^0 \stackrel{\partial_1}{\outof}
\Z\Lambda^0\oplus\Z\Lambda^0 \stackrel{\partial_2}{\outof}
\Z\Lambda^0 \outof 0$$
with $\partial_1= (1-M_1^t,1-M_2^t)$
and $\displaystyle \partial_2={ M_2^t-1 \choose 1-M_1^t }$ for a suitable choice of bases.

Convergence of the spectral sequence to $K_*(C^*(\Lambda))$ (Definition \ref{D:spectral-sequence}) and the above means that we have the following finite filtration of $\KK_1=K_1(C^*(\Lambda))$:
$$
0=F_0(\KK_1)\subseteq F_1(\KK_1)=F_2(\KK_1)=\KK_1,
$$
with $F_1(\KK_1)\cong H_1(\mathcal{D})$.  Hence, $K_1(C^*(\Lambda))\cong H_1(\mathcal{D})$ as required.

Now, we could proceed to obtain an isomorphism of $K_0(C^*(\Lambda))$ by use of the spectral sequence; however we choose to use the Pimsner-Voiculescu sequence in order to deduce relatively easily the action of the isomorphism on basis elements.

By applying the Pimsner-Voiculescu to $B\rtimes\Z$ and $B\rtimes\Z^2\cong (B\rtimes\Z)\rtimes\Z$ in succession, we may deduce that $K_0(i_B):K_0(B)\into K_0(B\rtimes\Z^2)$ factors through an injection $\Phi_1:H_0(\mathcal{B})\into K_0(B\rtimes\Z^2)$, where $i_B$ is the canonical injection $i_B:B\into  B\rtimes\Z^2$ and $\mathcal{B}$ is the chain complex defined in Lemma \ref{L:Koszul}.  We may also deduce that there is an exact sequence that constitutes the first row of the following commutative diagram:
$$\begin{CD}
0 @>>> H_0(\mathcal{B}) @>\Phi_1>> K_0(B\rtimes\Z^2) @>>> H_2(\mathcal{B}) @>>> 0\\
@. @V\Phi_2VV @VV\Phi_0V @VVV @. \\
0 @>\Phi>> H_0(\mathcal{D}) @>>> K_0(C^*(\Lambda)) @>>> H_2(\mathcal{D}) @>>> 0.
\end{CD}$$
where all downward arrows are isomorphisms.  In particular, $\Phi_0:K_0 (C^*(\Lambda)\into K_0(B\rtimes\Z^2)$ is the isomorphism constructed in Theorem \ref{thrm:K(C*(Lambda))-K(Crossed-Product)} and $\Phi_2:H_0(\mathcal{B})\into H_0(\mathcal{D})$ is one of the isomorphisms in Theorem \ref{T:Hom}.

Now $H_2(\mathcal{D})\subseteq\Z\Lambda^0$ is a free abelian group, thus the exact sequences split and we have an isomorphism $\Phi:H_0(\mathcal{D})\oplus H_2(\mathcal{D}) \into C^*(\Lambda)$ such that $\Phi(g\oplus 0)= \Phi_0\Phi_1\Phi_2^{-1}(g)$ for all $g\in H_0(\mathcal{D})$.  It is straightforward to check that $\Phi(\delta_u + \im\partial_1\oplus 0)=[p_u]$ for all $u\in\Lambda^0$.  Thus the Theorem is proved.
\end{proof}

Evidently complications arise when $k>2$, however it is worth noting
that under some extra assumptions on the vertex matrices it is
possible to determine a fair amount about the $K$-groups of higher rank graph
$C^*$-algebras. For example, the case $k=3$ is considered below.

\begin{prop}\label{P:k=3} Let $\Lambda$ be a row-finite 3-graph with
  no sources. Consider the following group homomorphisms defined by
  block matrices:
\begin{eqnarray*} \partial_1 &=& (1-M_1^t \;\; 1-M_2^t \;\; 1-M_3^t):\bigoplus_{i=1}^3\Z\Lambda^0 \into
  \Z\Lambda^0,\\ \partial_2 &=& \left(\!\!\!\begin{array}{ccc}M_2^t-1 &
      M_3^t - 1 & 0 \\ 1-M_1^t & 0 & M_3^t -1 \\ 0 & 1-M_1^t & 1- M_2^t\end{array}\!\!\!\right):\bigoplus_{i=1}^3\Z\Lambda^0 \into
  \bigoplus_{i=1}^3\Z\Lambda^0, \\ \partial_3 &=&
  \left(\!\!\!\begin{array}{c} 1-M_3^t \\ M_2^t -1 \\ 1-M_1^t
    \end{array}\!\!\!\right) : \Z\Lambda^0 \into \bigoplus_{i=1}^3 \Z\Lambda^0.
\end{eqnarray*}
There exists a short exact sequence:
$$
0 \into \coker \partial_1/G_0 \into K_0(C^*(\Lambda)) \into
\ker\partial_2/\im\partial_3\into 0,
$$
and
$$K_1(C^*(\Lambda))\cong \ker\partial_1/\im\partial_2 \oplus G_1,$$
where $G_0$ is a subgroup of $\coker\partial_1$ and $G_1$ is a subgroup of $\ker\partial_3$.
\end{prop}
\begin{proof}
By Theorem \ref{T:Main}, there exist short exact sequences
\begin{eqnarray*} 0 \into E^4_{0,0} \into K_0(C^*(\Lambda)) \into
  E^4_{2,-2}\into 0,\\ 0 \into E^4_{1,0} \into K_1(C^*(\Lambda)) \into
  E^4_{3,-2}\into 0. \end{eqnarray*}
However, since $E^4_{p,q}=0$ if $p\in\Z\backslash\{0,1,2,3\}$ the only
non-zero components of the differential $d^3$ are
$d^3_{3,q}:E^3_{3,q}\into E^3_{0,q+2}$, where $q\in 2\Z$. Moreover, as in the proof of Proposition \ref{P:k=2}, the
differential $d^2$ is the zero map.  Thus we
have
$$\begin{array}{ll}
E^4_{1,0}\cong E^3_{1,0}\cong E^2_{1,0}\cong \Homo_1(\mathcal{D}), &E^3_{0,0}\cong E^2_{0,0}\cong H_0(\mathcal{D}), \\
E^4_{2,-2}\cong E^3_{2,-2} \cong E^2_{2,-2}\cong\Homo_2(\mathcal{D}),& E^3_{3,-2}\cong E^2_{3,-2}\cong H_3(\mathcal{D}).
\end{array}
$$
Also note that $E^4_{3,-2}$ is isomorphic to $\ker d^3_{3,-2}\subseteq E^3_{3,-2}$, which is a subgroup of the free abelian group $H_3(\mathcal{D})\cong \ker \partial_3$.  Thus $E^4_{3,-2}$ is itself a free abelian group, from which we deduce that the exact sequence for
$K_1(C^*(\Lambda))$ splits.  Hence the result follows by setting $G_0$ to be the image of $\im d^3_{3,-2}\subseteq E^3_{0,0}$ under the isomorphism $E^3_{0,0}\into E^2_{0,0}\into H_0(\mathcal{D})$, and $G_1$ to be the image of $\ker d^3_{3,-2}\subseteq E^3_{3,-2}$ under the isomorphism $E^3_{3,-2}\into E^2_{3,-2}\into H_3(\mathcal{D})$.
\end{proof}
Now we consider two cases for which we can describe the $K$-groups of the $C^*$-algebra of a row finite $3$-graph with no sources in terms of its vertex matrices by the immediate application of Proposition \ref{P:k=3}.

\begin{cor}\label{C:k=3}
In addition to the hypothesis of Proposition \ref{P:k=3}:
\begin{enumerate}\renewcommand{\labelenumi}{(\arabic{enumi})}
\item if $\partial_1$ is
surjective then
\begin{eqnarray*} K_0(C^*(\Lambda)) &\cong& \ker \partial_2/ \im
  \partial_3, \\ K_1(C^*(\Lambda)) &\cong& \ker
  \partial_1/\im\partial_2 \oplus \ker \partial_3;
\end{eqnarray*}
\item if $\cap_{i=1}^3 \ker (1-M_i^t)=0$ then
\begin{eqnarray*} K_1(C^*(\Lambda)) &\cong&
  \ker \partial_1/\im \partial_2\end{eqnarray*} and there exists a
  short exact sequence
$$ 0 \into \coker \partial_1 \into K_0(C^*(\Lambda)) \into
\ker\partial_2/\im\partial_3\into 0. $$
\end{enumerate}
\end{cor}
\begin{proof}
To prove (1), note that we have $0=\coker\partial_1$ thus the exact sequence for $K_0(C^*(\Lambda))$ collapses to give the result for $K_0(C^*(\Lambda))$.  Also note that $0=\coker\partial_1=H_0(\mathcal{D})\cong E^3_{0,0}$.  Therefore $ d^3_{3,-2}:E^3_{3,-2}\into E^3_{3,0}$ is the zero map and $\ker d^3_{3,-2}=E^3_{3,-2}\cong E^2_{3,-2}\cong H_3(\mathcal{D})=\ker \partial_3$.  Therefore, $G_1$ in Proposition \ref{P:k=3} is $\ker\partial_3$ and (1) is proved.

To prove (2), if $\bigcap_{i=1}^3 \ker(1-M_i^t)=0$ then $\ker \partial_3=0$, which implies that $G_1$ in Proposition \ref{P:k=3} is the trivial group.  It also follows that $E^3_{3,-2}=0$ so that $\im d^3_{3,-2}=0$ and $G_0$ in Proposition \ref{P:k=3} is the trivial group.  Whence (2) follows immediately from Proposition \ref{P:k=3}.
\end{proof}
\begin{rem}${}$\\[-10pt]
\begin{enumerate}\renewcommand{\labelenumi}{(\roman{enumi})}
\item One may recover \cite[Theorem 3.1]{Pa97} from Theorem
  \ref{T:Main} by setting $k$ equal to 1.
\item By \cite[Corollary 3.5 (ii)]{KP00} a rank
$k$ Cuntz-Krieger algebra (\cite{RS99,RS01}) is isomorphic to a
$k$-graph $C^*$-algebra. Thus, Propsition \ref{P:k=2} generalises \cite[Proposition
  4.1]{RS01}, the proof of which inspired the methods used throughout
  this paper.
\item By showing that the $C^*$-algebra of a row-finite 2-graph,
  $\Lambda$, with
  no sources and finite vertex set, satisfying some further conditions, is isomorphic to a rank 2 Cuntz-Krieger algebra,
  Allen, Pask and Sims used Robertson and Steger's
  \cite[Proposition 4.1]{RS01}
  result to calculate the
  $K$-groups of $C^*(\Lambda)$ \cite[Theorem 4.1]{APS04}. Moreover, in
  \cite[Remark 4.7. (1)]{APS04}
  they note that their formulae for the $K$-groups holds for more
  general 2-graph $C^*$-algebras, namely the $C^*$-algebras of
  row-finite 2-graphs, $\Lambda$, with no sinks
  (i.e. $s^{-1}(v)\cap\Lambda^n\ne\emptyset$ for all
  $n\in\N^k,\;v\in\Lambda^0$) nor sources and finite vertex set.
  \item The notion of associating a $C^*$-algebra, $C^*(\Lambda)$, to a $k$-graph $\Lambda$ was generalised by Raeburn, Sims and Yeend \cite{RSY03} to include the case where $\Lambda$ is finitely-aligned; a property identified by them to enable an appropriate $C^*$-algebra to be constructed.  The family of finitely-aligned $k$-graphs and their associated $C^*$-algebras admit $k$-graphs with no sources and those that are not row-finite.  In \cite{F06}, Farthing devised a method of constructing, from an arbitrary finitely-aligned $k$-graph $\Lambda$ with sources, a row-finite $k$-graph with no sources, $\bar{\Lambda}$, which contains $\Lambda$ as a subgraph.  If, in addition, $\Lambda$ is row-finite then Farthing showed that $C^*(\bar{\Lambda})$ is strong Morita equivalent to $C^*(\Lambda)$ and thus has isomorphic $K$-groups to those of $C^*(\Lambda)$.  Therefore, in principal, the results in this paper could be extended to the case where $\Lambda$ is row-finite but with sources.
\end{enumerate}
\end{rem}

\section{The $K$-groups of unital $k$-graph $C^*$-algebras}\label{S:unital}

Recall that if $\Lambda$ is a row-finite higher rank graph with no
  sources then $\Lambda^0$ finite is equivalent to $C^*(\Lambda)$
  being unital (\cite[Remarks 1.6 (v)]{KP00}). Thus in this section we
  specialise in the case where the vertex set of our higher rank
  graph, hence each vertex matrix, is finite. We will continue to denote the Kasparov
spectral sequence converging to $K_*(C^*(\Lambda))$ of the previous section
by $\{(E^r,d^r)\}$ and we shall denote the torsion-free rank of an
abelian group $G$ by $r_0(G)$ (see e.g. \cite{Fu70}).

\begin{prop}\label{P:tor-free} If $\Lambda$ is a row-finite higher rank graph with no
  sources and $\Lambda^0$ finite then $K_0(C^*(\Lambda))$ and $K_1(C^*(\Lambda))$ have equal
  torsion-free rank.\end{prop}

\begin{proof} Let the rank of the given higher rank graph $\Lambda$ be $k$ and
let $|\Lambda^0|=n$.

Since,
$E^\infty_{p,q}\cong E^{k+1}_{p,q}$ for all $p,q\in\Z$ and $E^{k+1}_{p,q}=0$ if
$p\in\Z\backslash\{0,1,\ldots,k\}$ or $q$ odd by Lemma
\ref{L:E^infty}, it follows from the definition of convergence of
$\{(E^r,d^r)\}$ (Definition \ref{D:spectral-sequence}) that there exist finite filtrations,
$$
0=F_{-1}(\KK_0) \subseteq E^{k+1}_{0,0} \cong F_0(\KK_0)\subseteq F_1(\KK_0) \subseteq \cdots \subseteq
  F_{k-1}(\KK_0)\subseteq F_k(\KK_0)=\KK_0,
$$ and
$$
0=F_0(\KK_1) \subseteq  E^{k+1}_{1,0} \cong F_1(\KK_1)\subseteq
  F_2(\KK_1) \subseteq \cdots \subseteq F_{k-1}(\KK_1) \subseteq F_k(\KK_1)=\KK_1
$$
of $\KK_0=K_0(C^*(\Lambda))$ and $\KK_1=K_1(C^*(\Lambda))$ respectively,
such that $$E^{k+1}_{p,q}\cong F_p(\KK_{p+q})/F_{p-1}(\KK_{p+q}).$$

Thus,
\begin{eqnarray*}r_0(K_0(C^*(\Lambda))) &=&r_0(F_k(\KK_0)) =
r_0( F_{k-1}(\KK_0)) + r_0( E^{k+1}_{k,-k}) = \cdots \\ &=& r_0(F_0(\KK_0)) +
\sum_{s \ge 1} r_0(E^{k+1}_{s,-s}) = \sum_{s\in\Z}
r_0(E^{k+1}_{s,-s}),\end{eqnarray*} and
\begin{eqnarray*}r_0(K_1(C^*(\Lambda))) &=& r_0(F_k(\KK_1)) =
r_0( F_{k-1}(\KK_1)) + r_0( E^{k+1}_{k,-k+1}) = \cdots \\ &=& r_0(F_1(\KK_1)) +
\sum_{s \ge 2} r_0(E^{k+1}_{s,-s+1}) = \sum_{s\in\Z}
r_0(E^{k+1}_{s,-s+1}).\end{eqnarray*}

Now we claim that
$$ \sum_{s\in\Z} r_0(E^{k+1}_{s,-s}) - r_0(E^{k+1}_{s,-s+1}) = \sum_{s\in\Z}
r_0(E^2_{s,-s}) -r_0(E^2_{s,-s+1}).$$ To see that this holds it is sufficient to prove that for all $r\ge 2$ we have $$ \sum_{s\in\Z} r_0(E^{r+1}_{s,-s}) - r_0(E^{r+1}_{s,-s+1}) = \sum_{s\in\Z}
r_0(E^r_{s,-s}) - r_0(E^r_{s,-s+1}). $$

Recall that for all $r\ge 1,\;p,q\in\Z, \; E^{r+1}_{p,q}\cong
Z(E^r)_{p,q}/B(E^r)_{p,q}$ where $Z(E^r)_{p,q}=\ker d^r_{p,q}$ and
$B(E^r)_{p,q}=\im d^r_{p+r,q-r+1}$. Thus
\begin{eqnarray*} r_0(E^{r+1}_{p,q}) &=&
r_0(Z(E^r)_{p,q}) - r_0(B(E^r)_{p,q})\\ &=& r_0(Z(E^r)_{p,q}) -
r_0(E^r_{p+r,q-r+1}) + r_0(Z(E^r)_{p+r,q-r+1})\end{eqnarray*} for all $r\ge
  1,\;p,q\in\Z$. Moreover, it follows from the definition of the
Kasparov  spectral sequence that given any $r\ge 1$ and  $p,q,q'\in\Z$ with $q=q' \mod 2$
  there exist isomorphisms $\rho:E^r_{p,q}\into E^r_{p,q'},\;
  \sigma:E^r_{p-r,q+r-1}\into E^r_{p-r,q'+r-1}$ such that
  $d^r_{p,q'}\circ\rho=\sigma\circ d^r_{p,q}$. Therefore,
\begin{eqnarray*} \sum_{s\in\Z}r_0(E^{r+1}_{s,-s}) -
  r_0(E^{r+1}_{s,-s+1}) &=& \sum_{s\in\Z}\left\{ r_0(Z(E^r)_{s,-s}) -
r_0(Z(E^r)_{s+r,-s-r+2})\right. \\  &-& r_0(Z(E^r)_{s,-s+1}) +
r_0(Z(E^r)_{s+r,-s-r+1})\\ &+&\left. r_0(E^r_{s+r,-s-r+2}) -
r_0(E^r_{s+r,-s-r+1})\right\}\\ &=& \sum_{s\in\Z} r_0(E^r_{s,-s}) -
r_0(E^r_{s,-s+1}) \end{eqnarray*}
for all $r\ge 1$.
Combining the above gives $$
r_0(K_0(C^*(\Lambda))) - r_0(K_1(C^*(\Lambda))) = \sum_{s\in\Z}
r_0(E^2_{s,-s}) - r_0(E^2_{s,-s+1}).$$
Now, recall that for all $p\in\Z$ and $q\in 2\Z$, $E^2_{p,q}\cong
H_p(\Z^k,K_0(B))\cong \ker \partial_p/\im \partial_{p+1}$ by
Theorem \ref{T:Hom}. Therefore,
\begin{eqnarray*}
    && r_0(K_0(C^*(\Lambda))) - r_0(K_1(C^*(\Lambda))) = \sum_{s\in\Z} r_0(E^2_{2s,-2s}) - r_0(E^2_{2s+1,-2s}) \\
    &=& \sum_{s\in\Z} r_0(\ker \partial_{2s}) - r_0(\im \partial_{2s+1}) - r_0(\ker\partial_{2s+1}) + r_0(\im \partial_{2s+2}) \\
    &=& \sum_{s\in\Z} r_0(\ker \partial_{2s}) - r_0\left(\left(\bigoplus_{N_{2s+1}}\Z\Lambda^0\right)/\ker\partial_{2s+1}\right) - r_0(\ker\partial_{2s+1})\\
    &+& r_0\left(\left(\bigoplus_{N_{2s+2}}\Z\Lambda^0\right)/\ker\partial_{2s+2}\right) \\
    &=&  \sum_{s\in\Z} r_0(\ker \partial_{2s}) - {k \choose 2s+1}n + r_0(\ker \partial_{2s+1}) - r_0(\ker \partial_{2s+1}) \\
    &+&{k \choose 2s+2}n -  r_0(\ker \partial_{2s+2}) \\
    &=& \sum_{s\in\Z} \left\{{k \choose 2s} - {k \choose 2s-1} \right\}n \\
    &=& \sum_{s\in\Z} \left\{ {k-1 \choose 2s} + {k-1 \choose 2s-1} - {k-1 \choose 2s-1} -
    {k-1 \choose 2s-2} \right\} n = 0.
\end{eqnarray*}
\end{proof}

\begin{cor} If $\Lambda$ is a row-finite higher rank graph with no
  sources and $\Lambda^0$ is finite then there exists a non-negative
integer  $r$ such that for $i=0,1$, $$K_i(C^*(\Lambda))\cong \Z^r\oplus T_i$$
for some finite group $T_i$, where $\Z^0:=\{0\}$.
\end{cor}

\begin{proof} It is well-known that if $B$ is a finitely generated
subgroup of an abelian
group $A$ such that $A/B$ is also finitely generated then $A$ must be
finitely generated too \cite{Fu70}. Now, for all $p,q\in\Z$, $E^{k+1}_{p,q}$ is isomorphic
to a sub-quotient
of the finitely generated abelian group $E^2_{p,q}\cong
H_p(\Z^k,K_q(B))$, therefore $E^{k+1}_{p,-p}$ is also finitely
generated. Moreover, $E^{k+1}_{0,i}\cong F_0(K_i(C^*(\Lambda)))$ and
for $p\in\{1,2,\ldots,k\},\; E^{k+1}_{p,-p+i}\cong
F_p(K_i(C^*(\Lambda)))/F_{p-1}(K_i(C^*(\Lambda)))$, which implies that
$K_i(C^*(\Lambda))=F_k(K_i(C^*(\Lambda)))$ is finitely
generated. The result follows from Proposition \ref{P:tor-free} by noting that every finitely generated
abelian group $A$ is isomorphic to the direct sum of a finite
group with $\Z^r$, where $r=r_0(A)$ (see e.g. \cite[Theorem 15.5]{Fu70}).\end{proof}

\begin{rem}
Note that it is well-known that when $k=1$ we always have $T_1=0$ in the above, i.e. $K_1(C^*(\Lambda))$ is torsion-free.  However, for $k>1$, $K_1(C^*(\Lambda))$ may contain torsion elements.
\end{rem}

Formulae for the torsion-free rank and
torsion parts of the $K$-groups of unital $C^*$-algebras of row-finite $2$-graphs with no sources
can be given in terms of the vertex matrices (cf.
\cite[Proposition 4.13]{RS01}). This we do in Proposition \ref{P:Formulae}
below.

\begin{prop}[cf. {\cite[Proposition 4.13]{RS01}}]\label{P:Formulae} Let $\Lambda$ be a row-finite 2-graph with no sources
  and finite vertex set. Then
\begin{eqnarray*} r_0(K_0(C^*(\Lambda))) &=& r_0(K_1(C^*(\Lambda))\\
&=&  r_0(\coker (1-M_1^t, 1-M_2^t)) + r_0(\coker (1-M_1, 1-M_2)),\\
\tor (K_0(C^*(\Lambda))) &\cong& \tor (\coker (1-M_1^t,1-M_2^t)),\\
\tor (K_1(C^*(\Lambda))) &\cong& \tor (\coker (1-M_1, 1-M_2)).
\end{eqnarray*}
\end{prop}

\begin{proof} We have already seen in Proposition \ref{P:tor-free} that the
torsion-free rank of the $K_0$-group and $K_1$-group of a
$k$-graph are equal so it is sufficient to calculate the
torsion-free rank of $K_0(C^*(\Lambda))$. Let $n:=|\Lambda^0|$. By
Proposition \ref{P:k=2} we have
\begin{eqnarray*} r_0(K_0(C^*(\Lambda))) &=& r_0(\coker (1-M_1^t,
1-M_2^t)) + r_0\left(\ker {1-M_1^t \choose 1-M_2^t}\right) \\ &=&
r_0(\coker( 1-M_1^t, 1-M_2^t)) + n - r_0(\im (1-M_1, 1-M_2)) \\
&=& r_0(\coker(1-M_1^t, 1-M_2^t)) + r_0(\coker(1-M_1, 1-M_2)).
\end{eqnarray*}
Furthermore, the assertion about the torsion part of
$K_0(C^*(\Lambda))$ is obvious. The torsion part of
$K_1(C^*(\Lambda))$ is given by $$\tor(K_1(C^*(\Lambda))) \cong
\tor \left(\ker (1-M_1^t, 1-M_2^t)/\im{ M_2^t - 1 \choose 1-M_1^t}\right),$$
which is clearly isomorphic to $\tor(\coker {M_2^t-1 \choose
1-M_1^t})$. However, by reduction to Smith normal forms,
$\coker{M_2^t-1 \choose 1-M_1^t}$ is isomorphic to $\coker(1-M_1,
1-M_2)$.\end{proof}

\begin{rem}
We note that, in the case where $\Lambda$ is a row-finite $3$-graph with no sources and finite vertex set, and with $\partial_1,\partial_2$ defined as in Proposition \ref{P:k=3},   it is straightforward to show that if $\partial_1$ is surjective then
$$K_0(C^*(\Lambda)) \cong K_1(C^*(\Lambda))\cong \Z^m,$$
where $m:=r_0(\ker\partial_2) - |\Lambda^0| = r_0(\coker\partial_2) - |\Lambda^0|$ (with $\Z^0:=0$).
\end{rem}

\section{Applications and Examples}\label{S:examples}
We begin this section with two corollaries to the results in the preceding section, which facilitate the classification of the $C^*$-algebras of row-finite $2$-graphs with no sources.  We then end the paper with some simple illustrative examples.

\begin{cor}\label{C:unit} Let $\Lambda$ be a row-finite $2$-graph
  with no sources, finite vertex set and vertex matrices $M_1$ and $M_2$. Then there exists an isomorphism
$$
\Phi:\coker(1-M_1^t,1-M_2^t)\oplus\ker{ M_2^t-1 \choose 1-M_1^t }\into K_0(C^*(\Lambda))
$$
such that $\Phi(e+\im \partial_0)=[1]$, where $e(v)=1$ for all
   $v\in\Lambda^0$.
\end{cor}

\begin{proof} Follows immediately from Proposition \ref{P:k=2} and the fact that $\sum_{u\in\Lambda^0} p_u =1$.\end{proof}

\begin{rem}\label{R:Kirchberg-Phillips}
We note that the $C^*$-algebra of a row-finite $k$-graph $\Lambda$, with no sources, is separable, nuclear and satisfies the UCT \cite{RSc87}.  If in addition the $C^*$-algebra is simple and purely infinite we say that it is a Kirchberg algebra, and note that by the Kirchberg-Phillips classification theorem (\cite{K,P00}) it is classifiable by its $K$-theory (see \cite[Theorem 5.5]{KP00}).  We also note that conditions on the underlying $k$-graph have been identified, which determine whether the $C^*$-algebra is simple (\cite[Proposition 4.8]{KP00}, \cite[Theorem 3.2]{RSi07}) and purely infinite (\cite[Proposition 8.8]{S06}).
\end{rem}

\begin{cor}\label{C:2-graphs-with-same-vertex-mxs}
Let $\Lambda$ and $\Delta$ be two row-finite $2$-graphs with no sources.  Furthermore, suppose that $C^*(\Lambda)$ and $C^*(\Delta)$ are both simple and purely infinite, and that $\Lambda$ and $\Delta$ share the same vertex matrices.  Then $C^*(\Lambda)\cong C^*(\Delta)$.
\end{cor}

\begin{proof}
Let $\Lambda$ and $\Delta$ are two $2$-graphs satisfying the hypothesis.  Then, by Proposition \ref{P:k=2} their $K$-groups are isomorphic.

Suppose that the vertex set of $\Lambda$ (and hence that of $\Delta$) is infinite.  Then $C^*(\Lambda)$ and $C^*(\Delta)$ are both non-unital, and thus stable, Kirchberg algebras with isomorphic $K$-groups.  Thus by the Kirchberg-Phillips classification theorem $C^*(\Lambda)\cong C^*(\Delta)$.

In the case where the vertex set of $\Lambda$ (and hence that of $\Delta$) is finite, $C^*(\Lambda)$ and $C^*(\Delta)$ are both unital Kirchberg algebras with isomorphic $K$-groups.  Furthermore, by Corollary \ref{C:unit} we see that the isomorphism of $K$-groups maps the $K_0$-class of the unit of one of the $C^*$-algebras onto that of the other.  Therefore, by the Kirchberg-Phillips classification theorem we conclude that $C^*(\Lambda)\cong C^*(\Delta)$ and the result is proved.

\end{proof}

\begin{exs}${}$\\[-10pt]
\begin{enumerate}
\item Let $\Lambda$ be a row-finite $2$-graph with no sources.  Suppose that the vertex matrices of $\Lambda$ are both equal to $M$, say.  By Proposition \ref{P:k=2} and \cite[Theorem 3.1]{Pa97} we have:
\begin{eqnarray*}K_i(C^*(\Lambda))&\cong& \Z\Lambda^0/\im(1-M^t) \oplus \ker (1-M^t)\\ &\cong& K_0(C^*(E))\oplus K_1(C^*(E)),
\end{eqnarray*}
for $i=1,2$, where $E$ is the $1$-graph with vertex matrix $M$.  The isomorphism for $K_0(C^*(\Lambda))$ is immediately obvious and that for $K_1(C^*(\Lambda))$ is given by
$$
{x \choose y} + \im {M^t-1 \choose 1-M^t} \mapsto (x + \im (1-M^t)) \oplus (x+y).
$$

\item Fix non-zero $n_1,n_2\in\N^2$ and let $\Lambda$ be a $2$-graph with one vertex and vertex matrices $M_1=(n_1)$ and $M_2=(n_2)$ respectively. By Proposition \ref{P:k=2}, we have $K_0(C^*(\Lambda))\cong K_1(C^*(\Lambda) \cong \Z/g\Z$, where $g$ is the greatest common divisor of $n_1-1$ and $n_2-1$.

Note that we recover the $K$-groups of tensor products of Cuntz algebras \cite{C77} by letting $\Lambda$ be the product $2$-graph of two $1$-graphs each with one vertex and a finite number of edges.  We also note that tensor products of two Cuntz algebras are not the only examples of $C^*$-algebras of such $2$-graphs $\Lambda$ (cf. \cite[\S6]{KP00}).  However, by Corollary \ref{C:2-graphs-with-same-vertex-mxs}, they are, up to $*$-isomorphism, the only examples of Kirchberg algebras arising from row-finite $2$-graphs with one vertex.

\item For each positive integer $n$, let $\mbf{O}_n$ be the 1-graph with 1 vertex, $\star$, and $n$ edges (i.e. morphisms of degree 1), $\alpha_1,\alpha_2,\ldots,\alpha_n$.  Let $c:\mbf{O}_3 \times \mbf{O}_n \into \Z$ be the unique functor that satisfies $c(\alpha_i,\star)=\delta_{i,1}$ $(i=1,2,3)$ and $c(\star,\alpha_i)=1$ $(i=1,\ldots,n)$. Define $\Lambda$ to be the $2$-graph $\Z\times_c (\mbf{O}_3\times \mbf{O}_n)$.  Let $T_i:=1-M_i^t$, where $M_1$ and $M_2$ are the vertex matrices of $\Lambda$.  Then
\begin{eqnarray}\label{E:T1}
T_1 \delta_u &=& -\delta_u-\delta_{u+1} \mbox{ and}\\ \label{E:T2}
T_2 \delta_u &=& \delta_u - n\delta_{u+1},
\end{eqnarray}
where we identify $\Lambda^0$ with $\Z$.

Clearly, $\ker{-T_2\choose T_1}=0$.  Now consider $\coker(T_1,T_2)$ and for each $g\in\Z\Lambda^0$ let $[g]$ be the image of $g$ under the natural homomorphism $\Z\Lambda^0\into \coker(T_1,T_2)$. By (\ref{E:T1}) and (\ref{E:T2}) we have $(n+1) [\delta_u]=0$.  Therefore, $\coker(T_1,T_2)$ is a cyclic group, generated by $[\delta_0]$ say, whose order divides $n+1$.  We claim that $\rho[\delta_0]\ne 0$ for each $\rho=1,\ldots,n$.  Suppose the contrary, then we have
$$ \rho\delta_0 = T_1x + T_2y$$
for some $x,y\in\Z\Lambda^0$ and $\rho\in\{1,\ldots,n\}$.  Thus, for each $u\in\Lambda^0$ we have
\begin{eqnarray*}
 \rho\delta_0(u)&=& -x(u)-x(u-1)+y(u)-ny(u-1).\\
\end{eqnarray*}
Since $x,y\in\Z\Lambda^0$, there exists $N$ such that $x(u)=y(u)=0$ if $|u|>N$, which we assume, without loss of generality, to be greater than zero.  It follows that
\begin{eqnarray*}
y(-N)&=&x(-N),\\
y(u)&=& x(u)+(n+1)\sum_{j=0}^{N-1+u}n^{N-1+u-j}x(-N+j),\mbox{ if }-N+1\le u\le -1,\\
y(u) &=& x(u) + (n+1)\sum_{j=0}^{N-1+u}n^{N-1+u-j}x(-N+j) + \rho n^{u},\mbox{ if } u\ge 0.
\end{eqnarray*}
Setting $u=N+1$, we arrive at the contradiction $(n+1)|\rho n^{N+1}$.  Therefore, by Proposition \ref{P:k=2} $K_0(C^*(\Lambda))\cong \Z/(n+1)\Z$.

Now we turn our attention to $\ker(T_1,T_2)/\im{-T_2 \choose T_1}$.  Suppose that $x\oplus y \in \ker(T_1,T_2)$, then
there exists $N$ such that $x(u)=y(u)=0$ if $|u|>N$,
\begin{eqnarray}\label{E:y-in-terms-of-x-1} y(-N) &=& x(-N)\mbox{ and}\\\label{E:y-in-terms-of-x-2}
y(u)&=& x(u) + (n+1)\sum_{j=-N}^{u-1}n^{u-1-i}x(j)\mbox{ for all } u\ge -N+1.
\end{eqnarray}
Let $P:\Z\Lambda^0\oplus\Z\Lambda^0\into \Z\Lambda^0$ be the projection onto the second component, i.e. $P(x\oplus y)=y$.  From (\ref{E:y-in-terms-of-x-1}) and (\ref{E:y-in-terms-of-x-2}) we see that $P$ is injective on $\ker(T_1,T_2)$ and thus induces an isomorphism $\ker(T_1,T_2)/\im{-T_2\choose T_1} \cong P(\ker(T_1,T_2))/\im T_1$.  Moreover, $P(\ker(T_1,T_2))=\{y\in\Z\Lambda^0 \;|\; \sum_{j\in\Z} (-1)^jy(j)=0\}$.  Now given $y\in P(\ker(T_1,T_2))$, define $z:\Lambda^0\into \Z$ by
\begin{eqnarray*}
z(u)&=& 0 \mbox{ if } u<-N,\\
z(u) &=& \sum_{j=-N}^u (-1)^{u-j+1}y(j) \mbox{ if } u\ge -N.
\end{eqnarray*}
Then it is straightforward to show that $z\in\Z\Lambda^0$ and $T_1 z=y$, and thus $\ker(T_1,T_2)/\im{-T_2\choose T_1}$ is the trivial group. Therefore, by Proposition \ref{P:k=2} $K_1(C^*(\Lambda))=0$.

Note that $\Lambda$ satisfies the hypotheses of \cite[Proposition 4.8]{KP00} and \cite[Proposition 8.8]{S06} and thus $C^*(\Lambda)$ is a (stable) Kirchberg algebra.  It now follows from the Kirchberg-Phillips classification theorem that $C^*(\Lambda)$ is $*$-isomorphic to the stabilized Cuntz algebra $\mathcal{O}_{n+2}\otimes\K$.

\item Let $c:\mbf{O}_3\times\mbf{O_3}\times\mbf{O_3}\into \Z_2$ be the unique functor that satisfies\footnote{We extend the definition of the product of two higher rank graphs (Examples \ref{exs:k-graph_constructions}.2) to the product higher rank graph of three higher rank graphs in the natural way.  Note that if $\Lambda_i$ is a $k_i$-graph for $i=1,2,3$ then both $(\Lambda_1 \times \Lambda_2)\times\Lambda_3$ and $\Lambda_1 \times (\Lambda_2 \times \Lambda_3)$ are clearly pairwise isomorphic as $(k_1+k_2+k_3)$-graphs to $\Lambda_1\times\Lambda_2\times\Lambda_3$ .}
$$
\begin{array}{rcl@{,\hspace{1cm}}rcl@{,\hspace{1cm}}rcl}
c(\alpha_1,\star,\star)&=&0 & c(\star,\alpha_1,\star)&=&0 & c(\star,\star,\alpha_1)&=&1,\\
c(\alpha_2,\star,\star)&=&0 & c(\star,\alpha_2,\star)&=&0 & c(\star,\star,\alpha_2)&=&1,\\
c(\alpha_3,\star,\star)&=&1 & c(\star,\alpha_3,\star)&=&1 & c(\star,\star,\alpha_3)&=&1.
\end{array}
$$
Then the vertex matrices of $\Lambda:=\Z_2\times_c (\mbf{O}_3\times\mbf{O_3}\times\mbf{O_3})$ are
$$
M_1=M_2=\left(\begin{array}{cc} 1 & 2 \\ 2 & 1 \end{array}\right),\qquad
M_3=\left(\begin{array}{cc} 0 & 3 \\ 3 & 0 \end{array}\right).
$$
Following the notation in Proposition \ref{P:k=3}, for $i=1,2,3$, we have
$$\begin{array}{rcl@{\hspace{5mm}}rcl}
\partial_1 &=& \left(\begin{array}{rrrrrr}   -1&  -1&  -1&  -1&   1&  -3 \\
  -1&  -1&  -1&  -1&  -3&   1 \end{array}\right),& \partial_3 &=& \left(\begin{array}{rr}    1&  -3 \\
  -3&   1 \\
   1&   1 \\
   1&   1 \\
  -1&  -1 \\
  -1&  -1  \end{array}\right).\\
\partial_2 &=& \left(\begin{array}{rrrrrr}   1&   1&  -1&   3&   0&   0 \\
   1&   1&   3&  -1&   0&   0 \\
  -1&  -1&   0&   0&  -1&   3 \\
  -1&  -1&   0&   0&   3&  -1 \\
   0&   0&  -1&  -1&  -1&  -1 \\
   0&   0&  -1&  -1&  -1&  -1  \end{array}\right),
\end{array}
$$
To compute the $K$-groups of $C^*(\Lambda)$ we reduce the relevant matrices to their Smith normal forms (for a (not necessarily square) matrix $M$ we shall denote its Smith normal form by $S(M)$).  In particular,
$$U_1\partial_1V_1 = S(\partial_1)=S(\partial_3)^t=\left(\begin{array}{cccccc} 1&  0&  0&  0&  0&  0 \\
0&  4&  0&  0&  0&  0  \end{array}\right)$$ for some invertible matrices $U_1,\;V_1$.
Thus, $\coker \partial_1 \cong \Z/4\Z$, $\ker\partial_3=0$ and Corollary \ref{C:k=3} can be applied to deduce that there exists a short exact sequence
$$ 0 \into \Z/4\Z \into K_0(C^*(\Lambda)) \into \ker\partial_2/\im \partial_3 \into 0$$ and $K_1(C^*(\Lambda))\cong \ker \partial_1/\im \partial_2$.  Now
$$
U_2\partial_2 V_2 = S(\partial_2)= \left(\begin{array}{cccccc}  1&  0&  0&  0&  0&  0 \\
  0&  1&  0&  0&  0&  0 \\
  0&  0&  4&  0&  0&  0 \\
  0&  0&  0&  4&  0&  0 \\
  0&  0&  0&  0&  0&  0 \\
  0&  0&  0&  0&  0&  0  \end{array}\right),
  $$
for some invertible matrices $U_2,\;V_2$, and we see that there exists an isomorphism (induced by $V_2^{-1}$ or by $U_1$) of $\ker \partial_2/\im \partial_3$ onto $\Z/4\Z$.  Furthermore, we can now see that $\ker\partial_1/\im\partial_2$ is isomorphic to $\Z/4\Z\oplus\Z/4\Z$.  Hence, $K_0(C^*(\Lambda))$ is a group order 16 and $K_1(C^*(\Lambda))\cong \Z/4\Z \oplus\Z/4\Z$.\footnote{It is well-known that there are, up to isomorphism, 5 abelian groups of order 16.}

We note that it is, perhaps, slightly surprising that $C^*(\Lambda)$ shares, at least, one of its $K$-groups with that of $\mathcal{O}_5\otimes\mathcal{O}_5\otimes\mathcal{O}_5$, given that this relationship is not obvious at the level of $C^*$-algebras.

\item Let $c:\mbf{O}_2\times\mbf{O_3}\times\mbf{O_3}\into \Z_2$ be the unique functor that satisfies
$$
\begin{array}{rcl@{\hspace{1cm}}rcl@{,\hspace{1cm}}rcl}
c(\alpha_1,\star,\star)&=&0, & c(\star,\alpha_1,\star)&=&0 & c(\star,\star,\alpha_1)&=&1,\\
c(\alpha_2,\star,\star)&=&1, & c(\star,\alpha_2,\star)&=&1 & c(\star,\star,\alpha_2)&=&1,\\
 &&& c(\star,\alpha_3,\star)&=&1 & c(\star,\star,\alpha_3)&=&1.
\end{array}
$$
Then the vertex matrices of $\Lambda:=\Z_2\times_c (\mbf{O}_3\times\mbf{O_3}\times\mbf{O_3})$ are
$$
M_1=\left(\begin{array}{cc} 1 & 1 \\ 1 & 1 \end{array}\right),\qquad
M_2=\left(\begin{array}{cc} 1 & 2 \\ 2 & 1 \end{array}\right),\qquad
M_3=\left(\begin{array}{cc} 0 & 3 \\ 3 & 0 \end{array}\right).
$$
Now, for $i=1,2,3$, we have
$$
\begin{array}{lcr@{\hspace{5mm}}lcr} \partial_1 &=& \left(\begin{array}{rrrrrr}     0&  -1&   0&  -2&   1&  -3 \\
  -1&   0&  -2&   0&  -3&   1   \end{array}\right),&\partial_3 &=& \left(\begin{array}{rr}      1&  -3 \\
  -3&   1 \\
   0&   2 \\
   2&   0 \\
   0&  -1 \\
  -1&   0  \end{array}\right).\\
\partial_2 &=& \left(\begin{array}{rrrrrr}       0&   2&  -1&   3&   0&   0 \\
   2&   0&   3&  -1&   0&   0 \\
   0&  -1&   0&   0&  -1&   3 \\
  -1&   0&   0&   0&   3&  -1 \\
   0&   0&   0&  -1&   0&  -2 \\
   0&   0&  -1&   0&  -2&   0 \end{array}\right),
\end{array}
$$
As in the previous example we compute the Smith normal form of $\partial_1$ (and hence that of $\partial_3$) first, to determine whether Corollary \ref{C:k=3} is applicable.  We find that
$$S(\partial_1)=S(\partial_3)^t=\left(\begin{array}{cccccc} 1 & 0 & 0 & 0 & 0 & 0 \\ 0 & 1 & 0 & 0 & 0 & 0 \end{array}\right),$$
 thus $\coker\partial_1\cong \ker \partial_3 \cong 0$ and we may apply Corollary \ref{C:k=3} to deduce that
$$K_0(C^*(\Lambda))=\ker \partial_2/\im \partial_3\quad\mbox{and}\quad K_1(C^*(\Lambda))=\ker\partial_1/\im\partial_2.$$
Now,
$$ S(\partial_2)=\left(\begin{array}{cccccc}   1&  0&  0&  0&  0&  0 \\
  0&  1&  0&  0&  0&  0 \\
  0&  0&  1&  0&  0&  0 \\
  0&  0&  0&  1&  0&  0 \\
  0&  0&  0&  0&  0&  0 \\
  0&  0&  0&  0&  0&  0  \end{array}\right).$$
It follows that both $\ker\partial_2/\im\partial_3$ and $\ker\partial_1/\im\partial_2$ are trivial.  Thus, the $K$-groups of $C^*(\Lambda)$ are isomorphic to those of the Cuntz algebra $\mathcal{O}_2$.  Furthermore, It is clear that $\Lambda$ satisfies the hypotheses of \cite[Proposition 4.8]{KP00} and \cite[Proposition 8.8]{S06}, and therefore $C^*(\Lambda)$ is an unital Kirchberg algebra (cf. Remarks \ref{R:Kirchberg-Phillips}).  Applying the Kirchberg-Phillips classification theorem, we conclude that $C^*(\Lambda)\cong \mathcal{O}_2$.

\end{enumerate}
\end{exs}

\ifx\undefined\bysame
\newcommand{\bysame}{\leavevmode\hbox to3em{\hrulefill}\,}
\fi

\end{document}